\documentclass[11pt]{article}
\usepackage{}
\usepackage[total={6in, 9in}]{geometry}
\usepackage{amsfonts}
\usepackage{amsthm}
\usepackage{amssymb}
\usepackage{amsmath,amsfonts}
\usepackage{mathrsfs}
\usepackage{cases}
\usepackage{latexsym,bm}
\usepackage{indentfirst}
\usepackage{color}
\usepackage{ifpdf}
\usepackage{graphicx}
\usepackage{psfrag}
\usepackage{dsfont}
\usepackage[ruled,vlined]{algorithm2e}
\usepackage{enumerate}
\usepackage[
pdfauthor={},
pdftitle={},
pdfstartview=XYZ,
bookmarks=true,
colorlinks=true,
linkcolor=blue,
urlcolor=blue,
citecolor=blue,
bookmarks=true,
linktocpage=true,
hyperindex=true
]{hyperref}

\usepackage{pgf,tikz,pgfplots}
\usepackage{tikz-3dplot}
\usepackage{pifont}
\usepackage{refcount}

\newtheoremstyle{nonitalic}
{\topsep}   
{\topsep}   
{\normalfont}  
{}         
{\bfseries} 
{.}        
{ }        
{}         

\theoremstyle{nonitalic}

\newtheorem{THM}{\textbf{Theorem}}[section]

\newtheorem{DEF}[THM]{\textbf{Definition}}
\newtheorem{LEM}[THM]{\textbf{Lemma}}
\newtheorem{CLA}{\textbf{Claim}}[section]
\newtheorem{CON}[THM]{\textbf{Conjecture}}
\newtheorem{COR}[THM]{\textbf{Corollary}}

\newtheorem{REM}{\textbf{Remark}}

\newcommand{\iC}{\overset{\leftharpoonup }{C}}

\newcommand{\oC}{\overset{\rightharpoonup }{C}}

\newcommand{\pf}{\noindent\textbf{Proof}.\quad}

\newcommand{\CC}{\mathcal{C}}

\begin{document}
	\title{Hamiltonian cycles in tough $(P_4 \cup P_1)$-free graphs}
	\author{
		Songling Shan \\
		\medskip 
		Auburn University, Auburn, AL  36849\\
		\medskip 
		{\tt szs0398@auburn.edu}
	}
	
	\date{\today}
	\maketitle

\begin{abstract} 
In 1973, Chv\'atal conjectured that there exists a constant $t_0$  such that every  $t_0$-tough graph on at least three vertices  is Hamiltonian. This conjecture has inspired extensive research and has been verified for several special classes of graphs. Notably, Jung in 1978 proved  that every 1-tough $P_4$-free graph on at least three vertices is Hamiltonian. However, the problem remains challenging even when restricted to graphs with no induced $P_4\cup P_1$, the disjoint union of a path on four vertices and a one-vertex path. In 2013, Nikoghosyan   conjectured that every 1-tough $(P_4\cup P_1)$-free graph on at least three vertices is Hamiltonian.  Later in 2015, Broersma   remarked that ``this question seems to be very hard to answer, even if we impose a higher toughness." He instead posed the following question: ``Is the general conjecture of Chv\'atal's true for $(P_4\cup P_1)$-free graphs?" We  provide  a positive answer to Broersma's question by establishing that every $23$-tough $(P_4\cup P_1)$-free graph on at least three vertices is Hamiltonian.

\medskip

\end{abstract}

\section{Introduction}

We consider only simple  graphs. 
Let $G$ be a graph.
Denote by $V(G)$ and  $E(G)$ the vertex set and edge set of $G$,
respectively. Let $v\in V(G)$, $S\subseteq V(G)$, and $H\subseteq G$. 
Then  $N_G(v)$   denotes the set of neighbors
of $v$ in $G$, $d_G(v):=|N_G(v)|$ is the degree of $v$ in $G$, 
and $\delta(G):=\min\{d_G(v): v\in V(G)\}$ is the minimum degree of  $G$. 
Define  $N_G(v, S)=N_G(v)\cap S$,  $d_G(v,S)=|N_G(v, S)|$, 
 $N_G(S)=(\bigcup_{x\in S}N_G(x))\setminus S$, and $N_G(S,T)=N_G(S)\cap T$ for some $T\subseteq V(G)$.  
We write $N_G(v, H)$,  $d_G(v,H)$, and $N_G(H,T)$  respectively for $N_G(v, V(H))$,  $d_G(v, V(H))$, and  $N_G(V(H),T)$. 
We use $G[S]$ and $G-S$ to denote the subgraphs of $G$ induced by  $S$ and $V(G)\setminus S$, respectively. 
For notational simplicity we write $G-x$ for $G-\{x\}$.
Let $V_1,
V_2\subseteq V(G)$ be two disjoint vertex sets. Then $E_G(V_1,V_2)$ is the set
of edges in $G$  with one endvertex  in $V_1$ and the other endvertex  in $V_2$. 
For $u,v\in V(G)$, we write $u\sim v$ if $u$ and $v$ are adjacent in $G$, 
and we write $u\not\sim v$ otherwise.  Given two positive integers $p$ and $q$, and two sequences of vertices $u_1, \ldots, u_p$ and $v_1, \ldots, v_q$, 
we write $u_1, \ldots, u_p\sim v_1, \ldots, v_q$ if it holds that $u_i\sim v_j$ for each $i\in[1,p]$ and each $j\in [1,q]$. 
Given a graph $R$, we say that  $G$ is \emph{$R$-free} if $G$ does not contain $R$ as an induced subgraph. 
For an integer $k\ge 2$, we use $kR$ to denote the disjoint union of $k$ copies of $R$. When we say that $G$ is $(R_1\cup R_2)$-free,  we  take  $(R_1\cup R_2)$ as the vertex-disjoint union of 
two graphs $R_1$ and $R_2$. We use $P_n$ to denote a path on $n$ vertices.   For two integers $a$ and $b$, let $[a,b]=\{i\in \mathbb{Z}\,:\,   a\le i \le b\}$. 
Throughout this paper,  if not specified, 
we will assume $t$ to be a nonnegative real number.

 Let   $c(G)$ denote  the number of components of a graph $G$. 
Given a graph $G$, the  \emph{toughness}  of $G$,  denoted $\tau(G)$, is $\min\{|S|/c(G-S): S\subseteq V(G), c(G-S) \ge 2\}$  if $G$ is not
a complete graph, and is defined to be $\infty$ otherwise.  A graph    is called  \emph{$t$-tough}
if  its toughness is at least $t$. 
This concept was introduced by Chv\'atal~\cite{chvatal-tough-c} in 1973.
It is  easy to see that   every cycle is 1-tough and so every Hamiltonian graph is 1-tough. 
 Conversely,
Chv\'atal~\cite{chvatal-tough-c} proposed the following well-known conjecture. 

\begin{CON}[Chv\'{a}tal's Toughness Conjecture] \label{Con:toughness-conjecture} 
		There exists a constant $t_0$ such that every $t_0$-tough graph on at least three vertices is Hamiltonian.
	\end{CON}
	
Bauer, Broersma and Veldman~\cite{Tough-counterE} have constructed
$t$-tough graphs that are not Hamiltonian for all $t < \frac{9}{4}$, so
$t_0$ must be at least $\frac{9}{4}$ if Chv\'atal's Toughness Conjecture is true.
The conjecture has
been verified  for certain classes of graphs including 
planar graphs, claw-free graphs, co-comparability graphs, and
chordal graphs. For a more comprehensive list of graph classes for which the conjecture holds, see  the survey article by Bauer, Broersma, and Schmeichel~\cite{Bauer2006}  in 2006.  Some recently  established families of 
graphs for which the conjecture holds  include $2K_2$-free graphs~\cite{2k2-tough,1706.09029,OS2021}, 
and $R$-free 
graphs if $R$ is a 4-vertex linear forest~\cite{MR3557210} or  $R\in \{P_2\cup P_3, P_3\cup 2P_1,  P_2\cup 3P_1, P_2\cup kP_1\}$~\cite{S2021, P3-union-2P1,  HG2021, MR4455928,  MR4676626, MR4658427}, where  
$k\ge 4$ is an integer. In general, the conjecture is still wide open.


Among  the special classes of graphs for which Chv\'{a}tal's Toughness Conjecture was verified, notabely, Jung in 1978~\cite{P4-free-1-tough}  showed that 
every $1$-tough $P_4$-free graph on at least three vertices is Hamiltonian. 
However, the conjecture  remains challenging even when restricted to graphs with no induced $P_4\cup P_1$.  Nikoghosyan~\cite{Nikoghosyan} in 2013 conjectured that every 1-tough $(P_4\cup P_1)$-free graph on at least three vertices is Hamiltonian. In a 2015 survey~\cite{MR3443703}, Broersma  remarked that ``This question seems to be very hard to answer, even if we impose a higher toughness." He instead posed the following question: ``Is the general conjecture of Chv\'atal's true for $(P_4\cup P_1)$-free graphs?"  This same question was also asked by Li and Broersma in~\cite{MR3557210}. 
In this paper, we answer this question positively  by    establishing  the following result. 

\begin{THM}\label{thm:main-result}
	Every $23$-tough $(P_4\cup P_1)$-free graph on at least three vertices is Hamiltonian. 
\end{THM}

The toughness bound of $23$  in Theorem~\ref{thm:main-result} is likely not optimal. We choose this specific parameter primarily to facilitate the proof technique.  
The remainder of this paper is organized as follows. In the next section, we establish necessary preliminaries and lemmas. In the final section, we prove Theorem~\ref{thm:main-result}.

\section{Preliminaries and Lemmas}

Note that if $G$ is a $(P_4\cup P_1)$-free graph and $S$ is a cutset of $G$, then  each component of $G-S$ is $P_4$-free. 
Let $G$ be a $t$-tough $(P_4\cup P_1)$-free graph, where $t\ge 23$. 
Our main  strategy  for constructing a Hamiltonian cycle in $G$ is as follows (there is one case that needs a different approach). We first identify a set $S$ in $G$ such that $G-S$ is $P_4$-free and  each vertex of $S$
 has at least $\frac{n}{t+1}$ neighbors within $V(G)\setminus S$. We then proceed to find a cycle $C$ in $G$ that covers all vertices 
 of $G-S$. This cycle $C$ is constructed by utilizing  vertices from $S$ to 
link  together path segments covering the vertices of $G-S$.
 Lastly,  the remaining vertices of $S$ are iteratively ``inserted'' into $C$, leveraging their large number of  neighbors within $V(C)$, to ultimately obtain a Hamiltonian cycle for $G$.

To support this approach, we dedicate the first subsection to exploring the properties of $P_4$-free graphs.  
 In the second subsection, we demonstrate the existence of a cycle covering the vertices of $G-S$, given the aforementioned set $S$. 
 Finally, in the last subsection, we present the construction of a Hamiltonian cycle assuming the existence of a suitable set $S$ within $G$.

We start with  some definition  and  a  property about $(P_4\cup P_1)$-free graphs. 

Let $G$ be a graph and  $S\subseteq V(G)$. The graph $G$ is \emph{Hamiltonian-connected} if $G$ has a Hamiltonian $(u,v)$-path for any two distinct 
vertices $u,v$, and $G$  is \emph{Hamiltonian-connected with respect to $S$} if $G$ has a Hamiltonian $(u,v)$-path for any two distinct 
vertices $u,v$ such that $|\{u,v\}\cap S| \le 1$. 
Let  $x\in S$.  We say that $x$  is \emph{complete to}  a subgraph $H$ of $G-S$ if $N_G(x, H)=V(H)$,
and we say that $x$ is \emph{connected to} $H$ if  $N_G(x, H) \ne \emptyset$.  
Now  we say $S$ is complete to $H$ if every vertex of $S$ is complete to $H$. 
 If $S$ is a cutset of 
$G$, then an element $x\in S$ is called a \emph{minimal element} of $S$ if $x$ is contained in a minimal cutset 
of $G$ that is a subset of $S$. As any cutset contains a minimal cutset, every cutset in $G$ has a minimal element. 

\begin{LEM}\label{lem:P4-U-P1-cut}
Let $G$ be a $(P_4\cup P_1)$-free graph and $S$ be a  cutset of $G$. 
For $x\in S$ and $y\in N_G(x, G-S)$,  if $G-S$ has a vertex $z$ such that $z\not\sim x$, $z\not\sim y$, and $G-S$ has a component 
containing neither $y$ nor $z$ but a neighbor of $x$, 
then $x$ is complete to every component of  $G-S$  to which $x$ is connected in $G$ possibly except the one containing $z$.  
As a consequence, if $S$ is a minimal cutset of $G$, then  $x$  is complete to every component of  $G-S$    possibly except the one containing $z$. 
\end{LEM}

\pf Let $D_z$ be the component of $G-S$ that contains the vertex $z$.  We first show that $x$ is complete to all the components   of $G-S$ that contain neither $y$ nor $z$ but   contain a neighbor of $x$. 
Assume to the contrary that $G-S$ has a component $R$ with $V(R)\cap \{y,z\} =\emptyset$ such that $x$ 
has in $G$ a neighbor and a non-neighbor from $R$.  We choose vertices $w,w^* \in V(R)$ such that $ww^*\in E(R)$ and $x\sim w$ but $x\not\sim w^*$ ($w$ and $w^*$ exist by the connectedness of $R$). 
Then $yxww^*$ and $z$ form an induced $P_4\cup P_1$ in $G$, a contradiction. 
Thus $x$ is complete to all the components of $G-S$ that contain neither $y$ nor $z$ 
but   contain a neighbor of $x$. 

We next show that if $y\not\in V(D_z)$, then $x$ is also complete to the component of $G-S$ containing $y$. 
By the assumption, we know that $G-S$ has a component, say $R'$, 
containing neither $y$ nor $z$ but a neighbor of $x$. We let $y'\in N_G(x, R')$.  
The rest argument follows the same idea as above with $y'$  playing the role of $y$
and the component of $G-S$ that contains $y$ playing the role of $R$. 
\qed

%
%
\subsection{Properties of $P_4$-free graphs}

A path $P$ connecting two vertices $u$ and $v$ is called 
a {$(u,v)$-path}, and we write $uPv$ or $vPu$ in order to specify the two endvertices of 
$P$.  If $x$ and $y$ are two vertices on a path $P$, then $xPy$ is the subpath of $P$
with endvertices as $x$ and $y$. 
Let $uPv$ and $xQy$ be two paths. If $vx$ is an edge, 
we write $uPvxQy$ as
the concatenation of $P$ and $Q$ through the edge $vx$.
Let $P$ be a  $(u,v)$-path  in $G$ 
and $x\in V(G)\setminus V(P)$. If $P$ has an edge $yz$, where $y$ is in the middle of $u$ and $z$ along $P$,  such that $x\sim y, z$, then we say that the path  $uPyxzPv$ 
is obtained from $P$ by \emph{inserting} $x$ between $y$ and $z$.

The first two statements in the lemma below are consequences of $P_4$-freeness.

\begin{LEM}\label{lem:P4-cut}
	Let $G$ be a $P_4$-free graph and $S$ be  a cutset of $G$ such that each vertex of $S$ is connected in $G$
	to at least two distinct components of $G-S$.  Then 
	\begin{enumerate}[(1)]
		\item  For every $x\in S$ and every component $D$ of $G-S$,  if $x$ is  connected to $D$, then $x$ is complete to $D$. 
		\item  Let $S^*\subseteq S$ be a minimal cutset of $G$. Then every vertex of $S^*$  is complete to $G-S^*$. 
		\item If $G$ is not a  complete bipartite graph with $s(G) \ge 0$, then $G-S$ has 
		a component of order at least two or $S$ is not an independent set in $G$. 
	\end{enumerate}
\end{LEM}

	\pf  We only prove the last statement. 
	If each component of $G-S$ has order one and $S$ is an independent set in $G$, then by Lemma~\ref{lem:P4-cut}(2), $S$ is a minimal cutset of $G$.
	Thus every vertex of $S$ is connected in $G$ to every component of $G-S$.  However, this implies that $G$ is a complete bipartite graph as $c(G-S)=|V(G)\setminus S|$. 
	\qed 
	
Let $G$ be a graph. We call 
$$
s(G) =\max\{c(G-S)-|S|: S\subseteq V(G), c(G-S) \ge 2\} 
$$
the \emph{scattering number} of $G$ if $G$ is not complete; otherwise $s(G)=-\infty$. A set $S\subseteq V(G)$ with $c(G-S)-|S|=s(G)$ and $c(G-S) \ge 2$ 
is called a \emph{scattering set} of $G$. 
The first two results below were proved by Jung in 1978~\cite{P4-free-1-tough}. 
\begin{THM}[\cite{P4-free-1-tough}]\label{thm:P4-hamilton}
	Let $G$ be a $P_4$-free graph. Then  
	\begin{enumerate}[(1)]
		\item  $G$ has a Hamiltonian path if and only if $s(G) \le 1$, 
		\item $G$ is Hamiltonian if and only if $s(G) \le 0$ and $|V(G)| \ge 3$,  
		\item $G$ is Hamiltonian-connected if and only if $s(G)<0$. 
	\end{enumerate}
\end{THM}

\begin{THM}[\cite{P4-free-1-tough}]\label{thm:P4-path-cover}
	Let $G$ be a $P_4$-free graph, $S$ be a maximum scattering set of $G$, and $v_1, v_2\in V(G)$ be two distinct vertices.
	Then $V(G)$
	can be covered by  $\max\{1, s(G)\}$ disjoint paths such that in case $v_1\not\in S$ or $s(G) \le 0$, the vertex $v_1$ is an endvertex of 
	one of those paths; in case $s(G)<0$, the path is a $(v_1,v_2)$-path. 
\end{THM}

Theorem~\ref{thm:P4-path-cover}  was a claim in~\cite{P4-free-1-tough} and was used 
to prove Theorem~\ref{thm:P4-hamilton}. We will 
apply  Theorem~\ref{thm:P4-path-cover}  in proving Theorem~\ref{thm:P4-free-scattering-0}.  Before that, 
we need some properties about a maximal scattering set in a graph. 

\begin{LEM}\label{lem:scattering-set}
Let $G$ be a    graph with $s(G)\ge 0$ and $S\subseteq V(G)$ be a maximal scattering set of $G$. 
Then  the following statements hold. 
\begin{enumerate}[(1)]
	\item  For every nonempty proper subset $S_1$ of $S$,  vertices of $S_1$   are connected in $G$  to at least 
	$|S_1|+1$ components of $G-S$. 
	\item  We have $s(D) \le 0$
	for each component $D$ of $G-S$. 
	\item  Suppose further that $G$ is $P_4$-free.  If  $S^*\subseteq S$, 
	then $S\setminus S^*$ is a maximal scattering set of $G-S^*$. 
\end{enumerate} 
\end{LEM}

 \pf Note that $c(G-S^*-(S\setminus S^*))=c(G-S)$. 
 
   For (1), suppose to the contrary that there exists a proper subset $S_1$ of $S$ such that 
 vertices of  $S_1$   are  connected in total to at most  $|S_1|$ components of $G-S$. 
 Then we have $$c(G-(S\setminus S_1)) -|S\setminus S_1|\ge c(G-S)-|S_1|+1-|S\setminus S_1|=s(G)+1.$$ This gives a contradiction to 
 the fact that $S$ is a scattering set of $G$. 
 
 For (2), if there exists a component $D$ of $G-S$ such that $s(D) \ge 1$, then we let $T$ be a scattering set of $D$. 
 It follows by the definition that $c(D-T) =  |T|+s(D)$. Then 
 we have $$c(G-(S\cup T)) -|S\cup T|\ge c(G-S)+|T|-|S\cup T|=s(G).$$ This gives a contradiction to 
 the fact that $S$ is a maximal scattering set of $G$. 
 
 For (3),  the statement holds trivially if $S^*=\emptyset$. Thus we assume that $S^* \ne \emptyset$. If $S^*=S$, then by  Statement (2), the empty set is a maximal scattering set of $G$. 
 Hence, we also assume that $S^*\ne S$. Thus $S^*$ is a nonempty proper subset of $S$. 
 Suppose to the contrary that $S\setminus S^*$ is not a maximal scattering set of $G-S^*$. 
 Let $T$ be a maximal scattering set of $G-S^*$.  
 If $T\subseteq S\setminus S^*$ ($T$ is a proper subset as $T\ne S\setminus S^*$), then  
 $(S\setminus S^*)\setminus T \ne \emptyset$, and  as $S$ is a maximal scattering set of $G$, by Statement (1), vertices of $(S\setminus S^*)\setminus T$ 
 are connected in $G$ to at least $|(S\setminus S^*)\setminus T|+1$ components of $G-S$. 
 Thus 
 \begin{eqnarray*}
 c(G-S^*-T)-|T|  &\le& c(G-S)-(|(S\setminus S^*)\setminus T|+1)+1-|T|\\
 &=&c(G-S)-|S|+|S^*|\\
 &=&c(G-S^*-(S\setminus S^*))-|S\setminus S^*|. 
 \end{eqnarray*}
 This gives a contradiction to $T$ being a  maximal scattering set of $G-S^*$.  
 
 Thus $T\not \subseteq S\setminus S^*$, and so $T\cap (V(G)\setminus S) \ne \emptyset$. 
Let $D$ be a component of $G-S$ such that $T\cap V(D) \ne \emptyset$. 
Assume that  $V(D)\setminus T \ne \emptyset$.  
If there is a vertex of $S\setminus S^*$ that is connected in $G$ to $D$ but is not contained in $T$, 
then by Lemma~\ref{lem:P4-cut}(1),  all vertices of $V(D)\cap  T$ are connected in $G-S^*$ to only one component of $G-S^*-T$, a contradiction to  Lemma~\ref{lem:scattering-set}(1).  
If all vertices of $S\setminus S^*$ that are  connected in $G$ to $D$ are  contained in $T$,
then by Lemma~\ref{lem:scattering-set}(2), we know that all vertices of $V(D)\cap T$ are connected in $G-S^*$ 
to at most $|V(D)\cap T|$ components of $G-S^*-T$, a contradiction to  Lemma~\ref{lem:scattering-set}(1).  
Thus we must have  $V(D)\subseteq T$ for any component $D$ of $G-S$ for which $V(D)\cap T \ne \emptyset$. 
We assume that there are in total $k$ components of $G-S$ whose vertices are all contained in $T$, where $k\in [1,c(G-S)]$. 
Then we have 
\begin{eqnarray*}
	c(G-S^*-T)-|T|  &\le& c(G-S)-k -(|(S\setminus S^*)\setminus T|+1)+1-|T|\\
	&=&c(G-S)-|S|-k+|S^*|\\
	&=&c(G-S^*-(S\setminus S^*))-|S\setminus S^*|-k\\
	&<& c(G-S^*-(S\setminus S^*))-|S\setminus S^*|. 
\end{eqnarray*}
This  gives a contradiction to $T$ being a scattering set of $G-S^*$.  
 \qed 

Let $G$ be a $P_4$-free graph. 
Theorem~\ref{thm:P4-hamilton}(3) states that   $G$ is Hamiltonian-connected if $s(G)<0$. When $s(G)=0$ and $G$ is not a balanced complete bipartite graph, 
  we show  below that $G$ is  
  Hamiltonian-connected with respect to a maximal scattering set $S$ of $G$.

\begin{THM}\label{thm:P4-free-scattering-0}
	Let $G$ be a $P_4$-free graph with $s(G)=0$ such that $G$ is not a balanced complete bipartite graph, and let $S\subseteq V(G)$ be a maximal scattering set of $G$. Then  $G$ is Hamiltonian-connected with respect to $S$. 
\end{THM}

\pf The proof is by induction on $n:=|V(G)|$.  The smallest $P_4$-free graph satisfying the conditions is obtained 
from $K_4$ by removing an edge, say $xy$, and a maximal scattering set $S$ consists of the two vertices from $V(G)\setminus\{x,y\}$. 
It is then easy to check that $G$ has a Hamiltonian path connecting any two vertices $u, v$ of $G$ if $|\{u,v\}\cap S|  \le 1$.

Thus we assume that $n\ge 5$.  Let $u,v\in V(G)$ be any two distinct vertices such that $|\{u,v\}\cap S|  \le 1$.  
 We assume, without loss of generality, that $u\not\in S$. 
Let $x\in S$ be a minimal element of $S$.  In particular,  if a minimal  element 
of $S$ has in $G$ a neighbor from $S$, we choose $x$ to be such one.  
Let $G^*=G-x$. 
Then we have that $s(G^*)=1$ and  that $S^*:=S\setminus \{x\}$ is a maximal scattering set of $G^*$ by Lemma~\ref{lem:scattering-set}(3).  By  Lemma~\ref{lem:P4-cut}(2), $x$ is complete to $G-S$. 
By Theorem~\ref{thm:P4-path-cover}, $G^*$ has a Hamiltonian path  $P$ with $u$ as one of its endvertices.   Since 
$s(G^*)=1$, it follows that none of the endvertices of  $P$  is from $S^*$ and each component of $P-S^*$
is a Hamiltonian path of one and exactly one component of $G-S$. 
We consider 
two cases in constructing a Hamiltonian $(u,v)$-path $Q$ of $G$ based on $P$. 

Suppose first that  the  other endvertex of $P$ 
is $v$.  Then as $G$ is not a balanced complete bipartite graph,  
we have that  one component of $G-S$ has at least two vertices 
or  $x$ is adjacent in $G$ to a vertex from $S$ by Lemma~\ref{lem:P4-cut}(3).    In the former case, as all the vertices from one common component of $G^*-S^*$ are located consecutively with 
each other  
on $P$, we let 
$y$ and $z$ be two vertices of a component of $G-S$ that are consecutive on $P$. 
Then we can insert $x$ in between $y$ and $z$ in getting $Q$.
In the latter  case, we let $y\in S$ such that $xy\in E(G)$. 
Then as $s(G^*)=1$, any neighbor $z$ of  $y$ on $P$ belongs to $G-S$. 
Then we can insert $x$ between $y$ and $z$ in getting  $Q$.

Suppose  next that  the  other endvertex of $P$ 
is $w$ with $w\ne v$.     If $v=x$, then $Q=uPwx$ is a desired Hamiltonian path of $G$. 
Thus we assume that $v\ne x$. 
Recall that  $w\in V(G^*)\setminus S^*$. 
Then $v$ is an internal vertex of $P$. We let $v_1$  
be the neighbor of    $v$ in  the path $uPv$. 
If $v_1\in V(G^*)\setminus S^*$ or $v_1\in S^*$ and $x\sim v_1$, we let $Q=uPv_1 xwPv$. 
If $v_1\in S^*$ and $x\not\sim v_1$, then by Lemma~\ref{lem:P4-cut}(2),  $v_1$ is also a minimal element of $S$. 
Now we let $Q^*=uPv_1wPv$ and insert $x$ in $Q^*$ the same way as in the   case where $P$ is a $(u,v)$-path. 
\qed 
\subsection{A cycle covering vertices of $G-S$}

In this subsection, we demonstrate the existence of a cycle in a  $4.5$-tough $(P_4 \cup P_1)$-free graph $G$ that covers  all vertices of  $G-S$, where $S$ is a minimal cutset of $G$. 
Our approach proceeds in three stages: (1) Leveraging  the toughness condition, 
for each  component  $D$ of $G-S$,  we ``match''  to it  some number  (related to $s(D)$)
of vertices  $S_D$ from $N_G(V(D))\cap S$ (Lemma~\ref{lem:generalized-K1r-matching}); 
 (2) Applying  Theorems~\ref{thm:P4-hamilton}, \ref{thm:P4-path-cover}, and~\ref{thm:P4-free-scattering-0}, we decompose $G-S$ into path segments. Crucially, the endvertices of each path segment are strategically chosen to  adjacent to 
 a distinct vertices   from  $S_D$ (Lemmas~\ref{lem:path-partner-system1} and~\ref{lem:path-partner-system2}); and 
(3) Exploiting the $(P_4 \cup P_1)$-free structure of $G$, we interconnect these path segments via their associated $S$-vertices, ultimately constructing the desired cycle that covers all vertices of $G-S$ (Lemma~\ref{lem:cycle-covering-G-S}).

We again start with some general definitions.  Let $G$ be a graph. Two edges of $G$ are \emph{independent} 
if they do not share any endvertices. A \emph{matching}  $M$ in $G$ is a set of independent edges.  
A vertex is  \emph{$M$-saturated} or \emph{$M$-covered} if the vertex is an endvertex of an edge of $M$. 
Otherwise, the  vertex is    \emph{$M$-unsaturated} or \emph{$M$-uncovered}.   We ususally do not 
distinguish between $M$  and the subgraph of $G$ induced on $M$. 
An \emph{$M$-alternating path}
is a path in $G$ with edges alternating between edges of $M$ and edges of $E(G)\setminus M$.  
 A \emph{star-matching} in  $G$  is a set of vertex-disjoint copies of stars. The vertices of degree at least 2 in a star-matching are called the \emph{centers} of the star-matching. In particular, if every star in a star-matching is isomorphic to $K_{1,r}$, where $r\geq 1$ is an integer, we call the star-matching a $K_{1,r}$-\emph{matching}.  Thus a matching is a $K_{1,1}$-matching. 
 For a star-matching $M$, we denote by $V(M)$ the set of vertices covered by $M$. And if $x,y\in V(M)$ and $xy\in E(M)$, we say $x$ is a \emph{partner} of $y$.
Let  $\{S, T\}$ be a partition of $V(G)$. We use $G[S,T]$ to denote the bipartite subgraph of $G$
between $S$ and $T$.

Let $G$ be a graph, $S$ be a cutset  of $G$, and $D_1,D_2,\ldots,D_{\ell}$
be all the components of $G-S$, where $\ell\ge 2$ is an integer. 
For each $D_i$,  
we let  $S_i=N_G(D_i,S)$ and  $H_i=G[V(D_i), S_i]$.  Let $r\ge 1$ be an integer. 

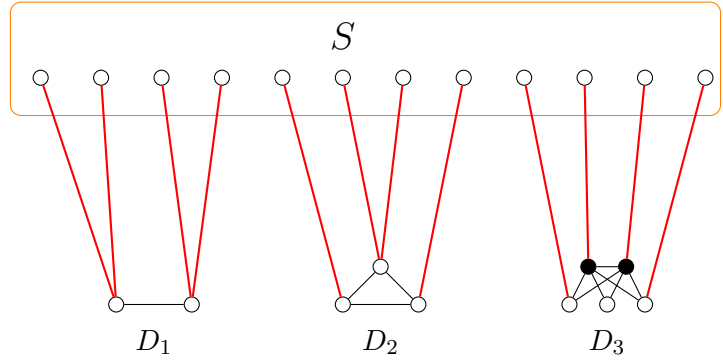
\begin{figure}[!htb]
\begin{center}
	\begin{tikzpicture}
		\draw[rounded corners, orange] (-4.4,4) rectangle (5,2.5);
		\node at (0,3.55) {\Large $S$};
		
		\foreach \i in {0,...,11} {
			\node[circle, draw, fill=white, inner sep=2pt] (s\i) at (-4+0.8*\i,3) {};
		}
		
		\node[circle, draw, fill=white, inner sep=2pt] (d1) at (-3,0) {};
		\node[circle, draw, fill=white, inner sep=2pt] (d2) at (-2,0) {};

		\draw (d1)--(d2);
		
		\node[circle, draw, fill=white, inner sep=2pt] (d4) at (0,0) {};
		\node[circle, draw, fill=white, inner sep=2pt] (d5) at (1,0) {};
		\node[circle, draw, fill=white, inner sep=2pt] (d6) at (0.5,0.5) {};
		
		\draw (d4) -- (d6) -- (d5) -- (d4);
		
		\node[circle, draw, fill=white, inner sep=2pt] (d7) at (3,0) {};
		\node[circle, draw, fill=white, inner sep=2pt] (d8) at (3.5,0) {};
		\node[circle, draw, fill=white, inner sep=2pt] (d9) at (4,0) {};
		\node[circle, draw, fill=black, inner sep=2pt] (d10) at (3.25,0.5) {};
		\node[circle, draw, fill=black, inner sep=2pt] (d11) at (3.75,0.5) {};
		
		\draw (d7) -- (d10) -- (d8) -- (d11) -- (d9) -- (d10);
		\draw (d7) -- (d11) -- (d10);
		
		\foreach \i in {0,1} {
			\draw[red, thick] (s\i) -- (d1);
		}
		\foreach \i in {2,3} {
			\draw [red, thick](s\i) -- (d2);
		}
		\draw [red, thick](s4)--(d4);
		\draw [red, thick](s5)--(d6)--(s6);
		\draw[red, thick] (s7)--(d5); 
		
	\draw[red, thick] (s8)--(d7);
	\draw[red, thick] (s9)--(d10);
	\draw[red, thick] (s10)--(d11);
	\draw [red, thick](s11)--(d9); 
		
		\node at (-2.5,-0.5) {$D_1$};
		\node at (0.5,-0.5) {$D_2$};
		\node at (3.5,-0.5) {$D_3$};
		
	\end{tikzpicture}
\end{center}
\caption{A depiction of a generalized $K_{1,4}$-matching, draw in red. In $D_3$, the set  $W$ consisting of 
the two black vertices is a cutset of $D_3$ such that $c(D_3-W)>|W|$.}
\label{f1}
\end{figure}

\begin{DEF}\label{def:good-matching}
For each bipartite graph  $H_i$, we let  $M_i$  be a star-matching of $H_i$.  Suppose $M_i$ satisfies the following properties: 
\begin{enumerate}[(M1)]
	\item $M_i$ has exactly $r$ edges;
	\item If $|V(D_i)| \ge r$, then $M_i$ is a matching; and if $|V(D_i)| <r$,   then  $M_i$ has exactly $|V(D_i)|$ components 
	such that each of the components is isomorphic to either  $K_{1, \lfloor r/|V(D_i)| \rfloor}$ or $K_{1, \lceil r/|V(D_i)| \rceil}$; 
	\item If   $D_i$  has a cutset $W_i$ such that $c(D_i-W_i) \ge |W_i|$,  then  $M_i$ covers 
	at least $\lfloor r/2 \rfloor$ vertices   from $V(D_i)\setminus W_i$.  Furthermore, if  $c(D_i-W_i) = |W_i|$, 
	each component of $D_i-W_i$ is  trivial, and $W_i$ is an independent set in $D_i$, 
	then   $M_i$  covers also a vertex of $W_i$. 
\end{enumerate}
Then we call $M_i$  a \emph{good star-matching} of $H_i$ with respect to $r$.

 For any $i, j\in [1,\ell]$, if there exists $S^*_i \subseteq S_i$  such that (i) $|S^*_i|=r$,
(ii) $S^*_i\cap S^*_j=\emptyset$ if $i\ne j$,  and (iii) $G[S^*_i,  V(D_i)]$ has  a good matching with respect to $r$,  then we say
that  $G$ has a \emph{generalized  $K_{1,r}$-matching}
with centers as  components of $G-S$, and call vertices in $S^*_i$ the \emph{partners}
of $D_i$ from $S$.   An example of a generalized $K_{1,4}$-matching is depicted in Figure~\ref{f1}. 
\end{DEF}

We will also need a theorem of K\"{o}nig on vertex covers. 
A \emph{vertex cover} in a graph is a set of vertices that contains an endvertex of every edge of the graph, 
and a vertex cover is \emph{minimum}  if its  size is minimum among  that of all vertex covers. 
The following classic result was due to K\"{o}nig.

\begin{THM}[\cite{Konig-thm}]\label{thm:Konig-Theorem}
	In any bipartite graph, the size of a maximum matching equals the size of  a minimum vertex cover. 
\end{THM}

Let $G$ be a graph,  $S\subseteq V(G)$, and  $D_1, \ldots, D_\ell$ be all 
the components of $G-S$ for some integer $\ell \ge 1$. For a rational number $t\ge 1$, we say that $G$
is $t$-tough with respect to   $S$ if for any cutset $W$ of $G$ for which 
$V(D_i)\setminus W\ne \emptyset$ for each  $i\in [1,\ell]$,  it holds that  $\frac{|W|}{c(G-W) }\ge t$. 
Note that $G$ is $t$-tough implies that $G$ is  $t$-tough with respect to   $S$ 
for any cutset $S$ of $G$.

\begin{LEM}\label{lem:generalized-K1r-matching} Let  $G$ be a graph, $t\ge 2$ be a rational number, and  $S$ be cutset  $G$. 
	If   $G$ is $t$-tough with respect to   $S$, then $G$ has a generalized  $K_{1,r}$-matching with centers as  components of $G-S$, where  $r=\lfloor t/2\rfloor $.  	
\end{LEM}

\pf As $S$ is a cutset of $G$, it is clear that every vertex of $V(G)\setminus S$ has in $G$ a non-neighbor.  Thus $G$ is $t$-tough with respect    $S$ implies that $d_G(v) \ge 2t$ for any $v\in V(G)\setminus S$.  
Let $D_1,D_2,\ldots,D_{\ell}$
be all the components of $G-S$, where $\ell \ge 2$ is an integer. 
For each $D_i$,  
we let  $S_i=N_G(D_i, S) $ and  $H_i=G[V(D_i), S_i]$. As  $G$ is $t$-tough with respect    $S$, we have $|S_i| \ge 2t$.

\begin{CLA}\label{claim:good-matching0}
	For each $i \in [1,\ell]$,  the bipartite graph   $H_i$ has a matching 
	of size  at least $\min\{|V(D_i)|, r\}$. 
\end{CLA}

\pf 
For otherwise, by Theorem~\ref{thm:Konig-Theorem}, a minimum vertex cover $Q$  of $H_i$ has size less than $\min\{|V(D_i)|, r\}$. 
Then $V(D_i)\setminus Q\ne \emptyset$, and as $|S_i|\ge 2t$, we know that $S\setminus Q\ne \emptyset$. 
However, $c(G-Q) \ge 2$ as there is no edge in $G$ between $D_i-Q$ and $G[S\setminus Q]$.   This gives 
a contradiction to  $G$ being $t$-tough with respect to $S$.
\qed 

\begin{CLA}\label{claim:good-matching}
	For each $i \in [1,\ell]$,  if    $H_i$ has a matching 
	of size  at least $\min\{|V(D_i)|, r\}$, then $H_i$  has a good star-matching with respect to $r$.
\end{CLA}

\pf   Let $M_i$ be a matching of  $H_i$ of size $\min\{|V(D_i)|, r\}$. 
  If $|V(D_i)| \ge r$, then $M_i$ satisfies (M1)-(M2) already. 
Thus we assume that $|V(D_i)| <r$ 
and so $|M_i|=|V(D_i)|$ by Claim~\ref{claim:good-matching0}. 
Then as $d_G(v) \ge 2t$ for every $v\in V(G)\setminus S$, we know that 
$d_G(v, S_i) \ge2t-|V(D_i)|> 2t-t/2>t/2$ for each $v\in V(D_i)$. Thus  for each $v\in V(D_i)$, 
we can choose a set  $T_v$ of $\lceil r/|V(D_i)|\rceil -1$   distinct vertices 
from $N_G(v,S_i\setminus V(M_i))$.  Furthermore, as $|N_G(v,S_i\setminus V(M_i))| >r$, 
for distinct $u, v\in V(D_i)$, we can choose $T_u$ and $T_v$ such that $T_u\cap T_v =\emptyset$. 
Then  $G[\left(V(M_i)\cap S_i\right) \cup (\bigcup_{v\in V(D_i)}T_v ),  V(D_i)]$ has a
star-matching that satisfies (M1)-(M2).


Next, we assume that $D_i$  has a cutset $W_i$ such that $c(D_i-W_i) \ge |W_i|$. 
It is clear that $|W_i| \le \frac{1}{2}|V(D_i)|$.   If $|V(D_i)| \le r$, then a star-matching of $H_i$ satisfying  properties (M1)-(M2)   also satisfies (M3).  Thus we assume that $|V(D_i)| >r$.  Thus  a star-matching  $M_i$ of $H_i$  satisfying  properties (M1)-(M2) is a matching of $H_i$.  We first show that $H_i$ has a matching covering at least $\lfloor \frac{r}{2} \rfloor$ vertices of $V(D_i)\setminus W_i$. 
If $|W_i| \le \lceil\frac{r}{2} \rceil $, then    $M_i$   is a desired matching already. 
Thus we assume that $|W_i| > \lceil\frac{r}{2} \rceil $.  
We show that $H_i^*=H_i[S_i, V(D_i)\setminus W_i]$
has a matching of size at least $\lfloor \frac{r}{2} \rfloor$.  For otherwise, by Theorem~\ref{thm:Konig-Theorem}, a minimum vertex cover $Q$  of $H_i^*$ has size less than  $\lfloor \frac{r}{2} \rfloor$. 
Then $(V(D_i)\setminus W_i)\setminus Q\ne \emptyset$, and as $|S_i|\ge 2t$, we know that $S\setminus Q\ne \emptyset$. 
However, $c(G-(Q\cup W_i)) \ge c(D_i-(Q\cup W_i))+1\ge |W_i|-|Q|+1 \ge 3$   as there is no edge in $G$ between $D_i- (Q\cup W_i)$ and $G[S\setminus Q]$.  
As 
\begin{equation}\label{eqn2}
\frac{|Q\cup W_i|}{c(D_i-(Q\cup W_i))} \le \frac{|Q\cup W_i|}{ |W_i|-|Q|}  =1+ \frac{2|Q|}{ |W_i|-|Q|} \le 1+ \frac{2(r-1)}{2} =r<t, 
\end{equation}
a contradiction to  $G$ being $t$-tough with respect to  $S$.  Thus $H_i^*$ has 
a matching  $M^*$ of size at least $\lfloor \frac{r}{2} \rfloor$. Since $H_i$ also has a matching  $M$  saturating at least $r$ vertices from $S_i$,     $H_i$ has a matching 
saturating all vertices in $V(M)\cap S_i$ and all vertices in $V(M^*) \cap V(D_i-W_i)$ (this follows by considering the symmetric difference of $M$ and $M^*$). 
Thus $H_i$  has a matching of size $r$ 
 that     covers at least $\frac{r}{2}$ vertices of $D_i-W_i$. 

 If  $c(D_i-W_i) = |W_i|$, 
each component of $D_i-W_i$ is  trivial, and $W_i$ is an independent set in $D_i$, then $V(D_i)\setminus W_i$ can also play the role of $W_i$. By the first part of (M3), 
we may assume that $M_i$ is a matching of $H_i$ of size $r$ that does not cover any vertex of $W_i$. 
Then by the same argument as above, we can find a matching $M^*$ of $H_i[S_i, W_i]$  of size  $\lfloor \frac{r}{2} \rfloor$. 
We  then  add edges of $M_i$ that are independent with edges of $M^*$  into $M^*$ to produce 
a size $r$ matching  of $H_i$ that   covers  $\lfloor \frac{r}{2} \rfloor$ vertices of $W_i$
and $\lceil \frac{r}{2} \rceil$ vertices of $V(D_i)\setminus W_i$ (as $M_i$  does not cover any vertex of $W_i$, it has at least $\lceil \frac{r}{2} \rceil$ 
edges that are independent with that of $M^*$). 

By the arguments above,  $H_i$ has a good star-matching with respect to $r$.  
\qed

\begin{CLA}\label{claim:Si-vertex-in-good-matching}
	For each $i \in [1,\ell]$, every vertex of $S_i$ is contained in a good star-matching  (with respect to $r$) of $H_i$. 
\end{CLA}

\pf  Let $M_i$ be a good star-matching (with respect to $r$) of $H_i$, and let $x\in S_i\setminus V(M_i)$. If $x$ is adjacent  in $G$ to a vertex $y\in V(M_i)\cap V(D_i)$,  then the star-matching obtained 
from $M_i$ by deleting an edge with one endvertex as $y$ and adding $xy$ is a star-matching  $M_i^*$ of size $r$ covering  $x$. 
It is clear that  $M_i$ is good with respect to $r$ implies that  $M_i^*$ is  also good with respect to $r$. 
If $x$ is adjacent  in $G$ to a vertex $y\in V(D_i)\setminus V(M_i)$, 
then we must have $|V(D_i)| >|M_i|$.  In case that $D_i$ has a cutset  $W_i$ such that $c(D_i-W_i) \ge |W_i|$, 
we choose an edge $uv\in M_i$ with $v\in S_i$  such that $u$ and $y$ 
are either both contained in $W_i$ or both contained in $V(D_i)\setminus W_i$.  Otherwise, we choose $uv\in M$ 
to be an arbitrary edge. 
Then the star-matching obtained 
from $M_i$ by deleting $uv$  and adding $xy$ is a good star-matching (with respect to $r$)  of $H_i$ covering $x$. 
\qed

By Claim~\ref{claim:Si-vertex-in-good-matching},  we let   $S_{i,1}, \ldots, S_{i,h_i}$, where $h_i \in \mathbb{N}$, be  all the  possible distinct subsets  of $S_i$ such that $|S_{i,j}|=r$, $\bigcup_{j=1}^{h_i} S_{i,j} =S_i$, and 
$G[V(D_i), S_{i,j}]$ has a good star-matching  with respect to $r$. 
Now we construct an $(r+1)$-uniform hypergraph $H$ based on $S$ and components of $G-S$. 
The hypergraph $H$ is bipartite with bipartition $S$ and  $\{d_1, \ldots, d_\ell\}$. 
For each $i\in [1,\ell]$ and the subsets $S_{i,1}, \ldots, S_{i,h_i}$ of $S_i$, we add $h_i$ hyperedges 
$S_{i,1}\cup \{d_i\}, \ldots, S_{i,h_i}\cup\{d_i\}$  containing $d_i$ to $H$.  

To finish the proof, it remains to show that  $H$ has a matching saturating 
$\{d_1, \ldots, d_\ell\}$.  Suppose not, we let $M$ be a maximum matching in $H$. Then $|M|\le \ell-1$. 
Without loss of generality, we let $d_1$ be an $M$-unsaturated vertex. Then by the same argument 
as  in the proof of  Hall's Theorem  on matchings in  bipartite graphs, we let 
$Z$ denote the set of all vertices connected to $d_1$ by $M$-alternating
paths. Since $M$ is a maximum matching, it follows 
that $d_1$ is the only $M$-unsaturated vertex in $Z$. Set $W=Z\cap \{d_1, \ldots, d_\ell\}$ and $T=Z\cap S$. 
Then  we have $|T|=r|W\setminus \{d_1\}|$ as there is a one-to-one correspondence given by $M$ 
between $W\setminus \{d_1\}$ and $|W|-1$ of  $r$-sets of $T$. 
Furthermore,  $H[W, S\setminus T]$ has no edge by $M$ being a maximum matching in $H$.

For any $d_i\in W$, by the maximality   of $M$, we know that $H[S_i\setminus V(M), V(D_i)]$ contains no edge. 
This implies that $G[S_i\setminus V(M), V(D_i)]$  has no good star-matching with respect to $r$.  
Then,  by Claim~\ref{claim:good-matching},  $G[S_i\setminus V(M), V(D_i)]$ has either 
no matching 
of size  $\min\{|V(D_i)|, r\}$, or it has  a matching of size   $\min\{|V(D_i)|, r\}$ 
but has no good-star matching with respect to $r$.   We define a subset $Q_i$ of $G[S_i\setminus V(M), V(D_i)]$ 
in three different cases below. 

If $G[S_i\setminus V(M), V(D_i)]$ has  
no matching 
of size at least $\min\{|V(D_i)|, r\}$, then by Theorem~\ref{thm:Konig-Theorem},  $H_i$ 
has a vertex cover $Q_i$  of size less than $\min\{|V(D_i)|, r\}$. 

Suppose now that  $G[S_i\setminus V(M), V(D_i)]$ has   a matching of size at least   $\min\{|V(D_i)|, r\}$ 
but has no good star-matching with respect to $r$.  
By the definition of a good star-matching, it follows that  $|V(D_i)|<r$ or $D_i$ has a 
cutset $W_i$ such that $c(D_i-W_i) \ge |W_i|$.  Let $M_i$ be a matching of $G[S_i\setminus V(M), V(D_i)]$ with size $\min\{|V(D_i)|, r\}$.

Assume   first  that $|V(D_i)| \ge r$ and  $D_i$ has a 
cutset $W_i$ such that $c(D_i-W_i) \ge |W_i|$.   By the same argument 
as in the proof of Claim~\ref{claim:good-matching},  we find 
a cutset $Q_i$ of 
$G[S_i\setminus V(M), V(D_i)]$ such that  $V(D_i)\setminus Q_i\ne \emptyset$ and $\frac{|Q_i|}{c(D_i-Q_i)} \le  r$ (see~\eqref{eqn2}).

Assume  then  that  $|V(D_i)|<r$.   Let $p$ 
be the principal remainder of $r$ divided by $|V(D_i)|$. 
For $p$  vertices $v\in V(D_i)$,  we let $F(v)$  be the set 
containing $\lceil r/|V(D_i)| \rceil$ 
duplications of $v$, and for the rest $|V(D_i)|-p$ 
vertices $v$ of $D_i$, we let $F(v)$  be the set 
containing $\lfloor  r/|V(D_i)| \rfloor$ 
duplications of $v$.  Let $T_i=\bigcup_{v\in V(D_i)} F(v)$. 
 We define     $H_i^*$ to 
be the bipartite graph  with bipartition $(S_i\setminus V(M), T_i)$, 
where   $e=xy$ with  $x\in S_i\setminus V(M)$ and $y\in F(v)$ for some $v\in V(D_i)$ 
is an edge of $H_i^*$ 
if and only if   $xv$ is an edge of $G[S_i\setminus V(M), V(D_i)]$. 
As there is no star-matching in 
$G[S_i\setminus V(M), V(D_i)]$ satisfying (M2),   it follows that $H_i^*$
has no matching of size $r$.  Then by Theorem~\ref{thm:Konig-Theorem},  $H_i^*$
has a  vertex cover $Q^*_i$  of size less than $r$. 
As all vertices from $F(v)$ for some $v\in V(D_i)$ have  the same 
neighbors in $H_i^*$ and $V(H_i^*)\setminus  Q_i^* \ne \emptyset$, 
it follows that $F(v)\cap Q_i^* = \emptyset$ for some $v\in V(D_i)$. 
Thus $G[S_i\setminus V(M), V(D_i)]$ has a subset $Q_i$
of  less than $r$ vertices such that $V(D_i)\setminus Q_i\ne \emptyset$ 
and there is no edge in $G$ between  $D_i-Q_i$
and $G[S_i\setminus (V(M) \cup Q_i)]$.

Assume, for notation convenience, that $W=\{d_1, \ldots, d_{|W|}\}$, 
and for some $k\in [1,|W|]$,  each of the components $D_1, \ldots, D_k$ 
has a cutset $Q_i$ defined as  in the first   case  right above (the case where $|V(D_i)| \ge r$ and  $D_i$ has a 
cutset $W_i$ such that $c(D_i-W_i) \ge |W_i|$).   Thus each $G[S_i\setminus V(M), V(D_i)]$ with $i\in [k+1, |W|]$ has a 
vertex cover  $Q_i$ with $|Q_i|<r$ such that $V(D_i)\setminus Q_i\ne \emptyset$. 
Let $q_i=c(D_i-Q_i)$ for each $i\in [1,k]$.  Then  we have $q_i\ge 2$ by~\eqref{eqn2}, and   
$|Q_i|\le rq_i$. 
Let $S^*=T \cup (\bigcup_{i=1}^{|W|}Q_i)$. 
Then we get 
\begin{eqnarray*}
	\frac{|S^*|}{c(G-S^*)} & \le&  \frac{|T|+ (r-1)(|W|-k)+rq_1+\ldots +rq_k}{|W|-k+q_1+\ldots+q_k} \\
	&\le & \frac{r(|W|-1)+ (r-1)(|W|-k)+ rq_1+\ldots +rq_k}{|W|+(q_1+\ldots+q_k-k)} \\
	&<& \frac{2r|W|+ 2rq_1+\ldots +2rq_k -r(q_1+\ldots q_k)}{|W|+(q_1+\ldots+q_k-k)} \\
	&\le & \frac{2r|W|+ 2r(q_1+\ldots +q_k -k)}{|W|+(q_1+\ldots+q_k-k)}   \le t, 
\end{eqnarray*}
giving a  contradiction to 
the   fact that  $G$ is $t$-tough with respect to   $S$. 
\qed

We will now construct paths that cover vertices of  some subgraph of a $(P_4\cup P_1)$-free graph. We need some basic definitions. 
\begin{DEF}
Let $G$ be a graph, $S\subseteq V(G)$,  $H\subseteq G-S$ be the union of some components of $G-S$.  Let $W=\emptyset$ if $s(H)\le 0$
and $W$ be a  maximal scattering set of $H$ otherwise.  

\begin{enumerate}[(1)]
	\item A \emph{path-cover}  $\mathcal{Q}$ of $H$  is the union of  some vertex-disjoint paths  
	such that $V(H)\subseteq V(\mathcal{Q})$.   
	\item A path-cover  $\mathcal{Q}$ of $H$  with components $R_1, \ldots, R_k (k\in \mathbb{Z})$  is a \emph{basic path-cover} of $H$ 
	if   $\mathcal{Q}$ satisfies the following conditions:
	\begin{itemize}
		\item $V(\mathcal{Q})=V(H)$, 
		\item $k =\max\{1, s(H)\}$, 
		\item $V(R_1)$ consists of  all vertices of $W$ 
		and  vertices of  $|W|+1$ components of $H-W$ (if $s(H) \ge 1$, this condition implies that all vertices from 
		the same component of $G-S$ form a subpath of $R_1$, and vertices of $W$ are used internally to link
		these $|W|+1$ subpaths), 
		\item  $H[V(R_i)]$ for each $i\in [2,k]$ is a 
		component of  $H-W$. 
	\end{itemize}
	\item A path-cover  $\mathcal{Q}$  of $H$ 
	is \emph{$S$-matched}  if the two endvertices of each path of $\mathcal{Q}$  belong  to $S$.   An \emph{$S$-vertex} of $\mathcal{Q}$
	is a vertex belonging to $V(\mathcal{Q})\cap S$, and  an \emph{$S$-endvertex}   is an $S$-vertex that is an endvertex of a component of $\mathcal{Q}$. 
	\item An $S$-matched path-cover  $\mathcal{Q}$  of $H$ is an \emph{$S$-matched basic path-cover}  if no two $S$-vertices are adjacent 
	in $\mathcal{Q}$ and $\mathcal{Q}-S$ is  a basic path-cover of $H$. 
	
	\item Let $\mathcal{Q}$  be an $S$-matched path-cover  of $H$.   
	Then  two components  $x_1u_1R_1v_1y_1$  and $x_2u_2R_2v_2y_2$ of $\mathcal{Q}$
	are \emph{linkable}  if   there exists $z\in \{u_2, v_2\}$, say $z=u_2$ such that 
	[($y_1\sim u_2$ or $x_2\sim v_1$) and ($y_2\sim u_1$ or $x_1\sim v_2$)] or 
	[($x_1\sim u_2$ or $x_2\sim u_1$) and ($y_2\sim v_1$ or $y_1\sim v_2$)]. 
	If $x_1u_1R_1v_1y_1$  and $x_2u_2R_2v_2y_2$ are linkable,  say $y_1\sim u_2$ and $y_2\sim u_1$, 
	then the new path  $x_1u_1R_1v_1y_1u_2R_2v_2y_2$ is called 
	a \emph{link} of $x_1u_1R_1v_1y_1$  and $x_2u_2R_2v_2y_2$. 
	
	\item Let  $\mathcal{Q}$ be  an $S$-matched  basic path-cover of $H$.    
	Then the  \emph{partner}  of an $S$-endvertex   is the neighbor of the $S$-vertex   in $\mathcal{Q}$. 
\end{enumerate}
\end{DEF}

By the definition of a basic path-cover, we have the following fact. 
\begin{REM}
	Let $\mathcal{Q}$  be an $S$-matched basic path-cover  of $H$ with $c(\mathcal{Q}) \ge 2$.   
	Then for any two components  $uPv$ and $xQy$ of $\mathcal{Q}$,  we have $E_G(N_P(\{u,v\}), N_Q(\{x,y\})) =\emptyset$ as 
	the vertices of $N_P(\{u,v\})$ and the vertices of  $N_Q(\{x,y\})$ are respectively  from two distinct components of $H$. 
	\end{REM}

Let $uPv$ and $xQy$ be two vertex-disjoint  paths and $z$ be a vertex not on $P$ or $Q$ such that $z\sim v, x$. 
We say that  \emph{linking  $P$ and $Q$ using $z$ in the order of $uPv$,  $xQy$}
consists of adding the edges  $zv$ and $zx$ to $P\cup Q$, 
 thereby obtaining the new path  $uPvzxQy$.

\begin{LEM}\label{lem:inserting}
	Let $G$ be a  $(P_4\cup P_1)$-free graph,   $S\subseteq V(G)$ such that $G-S$ is $P_4$-free,  and $D$ be a component of $G-S$. 
Suppose that  $s(D) \ge 0$ and $D$ is not a  complete bipartite graph. Let $W$ be 
a maximal scattering set of $D$, and $z\in W$ be a minimal element of $W$. 
Then if $\mathcal{Q}$ is an $S$-matched basic path-cover of $D-z$,  we can get 
an $S$-matched 
 path-cover of $D$ by either linking two components of $\mathcal{Q}-S$ 
using $z$  if $c(\mathcal{Q}) \ge 2$ or inserting $z$ into the component of $\mathcal{Q}$ if $c(\mathcal{Q}) =1$. 
\end{LEM}

 \pf   By Lemma~\ref{lem:scattering-set}(3), we have $s(D-z) \ge 1$. 
 Let $k=s(D-z)$, and $Q_1, \ldots, Q_k$ be all the components of $\mathcal{Q}$,
 where $Q_i=x_iu_iQ_iv_iy_i$ with $x_i, y_i \in S$, and $u_i, v_i\in V(D)$. 
 
 If $c(\mathcal{Q}) \ge 2$, then   $z \sim u_i,v_i$ for each $i\in[1,k]$ by Lemma~\ref{lem:P4-cut}(2).  
 Now  $$x_1u_1Q_1v_1z u_2Q_2v_2y_2, Q_3, \ldots, Q_k$$ form an $S$-matched basic path-cover of $D$. 
 
 If $c(\mathcal{Q}) =1$, then we have  $s(D)=0$ by Lemma~\ref{lem:scattering-set}(3). As $s(D-z) =1$, no two vertices of $W\setminus \{z\}$ are 
 consecutive on $Q_1$, and all the vertices from the same component of $D-z-W$ 
 are consecutive on $Q_1$. 
Since  $D$ is not a balanced complete bipartite graph,  $D-W$ 
has a component of order at least 2 or $D[W]$ has an edge by Lemma~\ref{lem:P4-cut}(3). 
In the former case,  we  insert $z$ on $Q_1$  in between two vertices of $D-W$ that 
are from the same component of $D-W$.  The resulting path is an $S$-matched  path-cover of $D$. 
 In the later case, we let $z_1z_2\in E(D[W])$.  
If $z$ is one of $z_1$ and $z_2$, say $z=z_1$, then 
we can insert $z_1$ between $z_2$ and one neighbor of $z_2$ on $Q_1$. 
The resulting path is an $S$-matched  path-cover of $D$. 
Thus we  assume that $z\not\in \{z_1,z_2\}$. Since $D$ is $P_4$-free 
and $z_1z_2\in E(D)$, if we let $\CC(z_i)$ be the set of components of $G-S$ 
that $z_i$ is connected to for each $i\in [1,2]$, then we must have $\CC(z_1) \subseteq \CC(z_2)$
or $\CC(z_2) \subseteq \CC(z_1)$. Without loss of generality, we assume $\CC(z_2) \subseteq \CC(z_1)$. 
We first replace $z_1$ by $z$ on $Q_1$, that is, deleting $z_1$ but joining $z$ to the two neighbors of $z_1$ 
on $Q_1$ to get $Q_1^*$, then we insert $z_1$ between $z_2$ and a neighbor of $z_2$ on $Q_1^*$. 
The resulting path is an $S$-matched  path-cover of $D$.
 \qed

\begin{LEM}\label{lem:path-partner-system1}
	Let $G$ be a   $(P_4 \cup P_1)$-free graph, and let $S\subseteq V(G)$ such that $G-S$ is $P_4$-free.  
	Then the following statements hold. 
\begin{enumerate}[(1)]
	\item  Suppose that $G$ is 4-tough with respect to $S$.  If $s(G-S) \ge 1$, then $G-S$ has an $S$-matched basic path-cover with $s(G-S)$ components. 
	\item Suppose that $G$ is 4.5-tough with respect to $S$.   If $s(G-S) \le 0$, then $G-S$ has an $S$-matched  path-cover with a single  component.  
\end{enumerate} 
\end{LEM}

\pf    For (1),  if $c(G-S)=1$, we let $D_1=G-S$, and let $S_1\subseteq V(D_1)$ be a maximal scattering set of $D_1$ and $\ell=1$. 
 If $c(G-S) \ge 2$, we 
let $D_1, \ldots, D_\ell$ be all the components of $G-S$, where $\ell:=c(G-S)$.  For each $D_i$,   let $S_i\subseteq V(D_i)$ be a maximal scattering set of $D_i$ if $s(D_i) \ge 1$, and let $S_i=\emptyset$ otherwise.  
Let $W=\bigcup_{i=1}^\ell S_i$. 
We apply induction on  $|W|$ in completing the proof. 

If $|W|=0$, then as $s(G-S) \ge 1$, the definition of $W$ and the condition that $s(G-S) \ge 1$  implies that $c(G-S) \ge 2$. 
Applying  Lemma~\ref{lem:generalized-K1r-matching},  we find a  generalized $K_{1,2}$-matching of $G$ with centers as components 
$D_1, \ldots, D_\ell$ of $G-S$.  In particular, each $D_i$ has two distinct partners $x_i, y_i$ from $S$  such that when $|V(D_i)|\ge 2$, there exist distinct $u_i, v_i\in V(D_i)$
for which $x_iu_i, y_iv_i\in E(G)$, and 
 $G[V(D_i), \{x_i,y_i\}]$ has a good star-matching with respect to $2$. 
For notation uniformity, when $D_i$ is a trivial component of $G-S$, we let $u_i=v_i$  be the vertex in $V(D_i)$. 
As $s(D_i) \le 0$ by  the assumption that $W=\emptyset$, 
 each $D_i$ is either Hamiltonian-connected,  a  balanced complete bipartite graph,  or Hamiltonian-connected 
 with respect to a cutset $W_i$ of $D_i$.  Since $G[V(D_i), \{x_i,y_i\}]$ has a good star-matching $\{x_iu_i, y_iv_i\}$, 
$D_i$  has a Hamiltonian $(u_i, v_i)$-path $P_i$.  Thus we get a path  $Q_i=x_iu_iP_iv_iy_i$, and so $Q_1, \ldots, Q_\ell$ is an $S$-matched basic path-cover of $G-S$.

Thus we assume that $|W| \ge 1$.  Without loss of generality, we assume that $S_1\ne \emptyset$.  This implies that $s(D_1) \ge 1$. 
Let $S_{11} \subseteq S_1$ be a minimal cutset of $D_1$.  Then we know that $D_1[S_{11}, V(D_1)\setminus S_{11}]$ is 
a complete bipartite graph, and $S\cup S_{11}$ is a cutset of $G$. 
Note that $S_1\setminus S_{11}$ is a maximal scattering set of $D_1-S_{11}$ by Lemma~\ref{lem:scattering-set}(3)
and $|W\setminus S_{11}|<|W|$.  
By induction,  $G-(S\cup S_{11})$ has an $(S\cup S_{11})$-matched basic  path-cover  $\mathcal{Q}$ with $s(G-(S\cup S_{11}))$ components.
In particular,  there are $s(D_1)+|S_{11}|$  components of $\mathcal{Q}$  that are covering vertices of $D_1-S_{11}$. 
We assume that these paths are  $Q_1:=x_1u_1R_1v_1y_1, \ldots, Q_k:=x_ku_kR_kv_ky_k$, where $k=s(D_1)+|S_{11}| \ge 1+|S_{11}|$,  $R_i:=u_iQ_iv_i$,  $x_i,y_i\in S\cup S_{11}$,   and $u_1R_1v_1$ 
is the path  containing   vertices of $S_1\setminus S_{11}$.  Among all these $k$ paths $Q_1, \ldots, Q_k$, at most $|S_{11}|$ of them 
that each contain a vertex of $S_{11}$.  As the endvertices of each $R_i$ are from $V(D_1)\setminus S_1$, and $D_1[S_{11}, V(D_1)\setminus S_{11}]$ is 
a complete bipartite graph, we know  each vertex of $S_{11}$  is adjacent in $G$ to all the endvertices of the paths $R_1, \ldots, R_k$.    We take $|S_{11}|$ paths from $Q_2, \ldots, Q_k$ such that 
all the paths that contain a vertex of $S_{11}$ are selected.   Without loss of generality, 
we let those paths be $Q_2, \ldots, Q_{p+1}$, where $p=|S_{11}|$.  As each path is matched to two vertices of $S\cup S_{11}$, 
there are two paths among $Q_1, \ldots, Q_{p+1}$ such that each of them has a partner from $S$. 
Let $Q_i$ and $Q_j$  be two paths with $i,j\in [1,p+1]$  and $i<j$ such that  one  vertex from $\{x_i, y_i\}$  and  one vertex from $\{x_j,y_j\}$  are in 
$S$. By  exchanging the labels of $x_i$ and $y_i$, and of $x_j$ and $y_j$ if necessary, we assume that 
$x_i, y_j\in S$.  
Then we link $R_1, \ldots, Q_i-y_i, \ldots, Q_j-x_j, \ldots,  R_{p+1}$  into one path $Q_1^*$ 
in the order of 
\begin{eqnarray*}
x_iQ_iv_i, u_1R_1v_1, \ldots, u_{i-1}R_{i-1}v_{i-1}, u_{i+1}R_{i+1}v_{i+1}, \ldots, \\ 
u_{j-1}R_{j-1}v_{j-1}, u_{j+1}R_{j+1}v_{j+1}, \ldots, u_{p+1}R_{p+1}v_{p+1}, u_jQ_jy_j
\end{eqnarray*}
by using  vertices  of $S_{11}$. 
Then $Q_1^*$ and the rest intact components of $\mathcal{Q}$ form  an $S$-matched basic path-cover of $G-S$.

For (2),  we let $x_1,y_1 \in S$ be distinct such that  $x_1,y_1$ are both adjacent to the unique vertex 
of $G-S$ if $|V(G-S)| =1$,  and  $x_1 \sim u_1$ and $y_1\sim v_1$ for distinct $u_1,v_1\in V(G-S)$  if $|V(G-S)|  \ge 2$. Such vertices $x_1,y_1, u_1,v_1$ exist as $G$ is 4-tough with respect to $S$. 

If $s(G-S) \le -1$ (this includes the case when $|V(G-S)|\in [1,2]$),  then $G-S$ has a Hamiltonian $(u_1,v_1)$-path $R_1$ by Theorem~\ref{thm:P4-hamilton}(3). Let $Q_1=x_1u_1R_1v_1y_1$. Then $Q_1$ is an $S$-matched basic path-cover with  one component. 

If $s(G-S) =0$ and $G-S$ is not a complete bipartite graph, 
we let $W\subseteq V(G-S)$ be a maximal scattering set of $G-S$, and let $w\in W$ be 
a minimal element of $W$. Then $s(G-w-S)=1$ by Lemma~\ref{lem:scattering-set}(3).  Since $G$ is 4.5-tough with respect to $S$, 
it follows that $G-w-S$ is 4-tough with respect to $S$. 
Applying Lemma~\ref{lem:path-partner-system1}(1), 
$G-w-S$ has an $S$-matched basic path-cover with  one component.  By  Lemma~\ref{lem:inserting}, $G-S$ has an $S$-matched  path-cover with  one component.

Thus we assume that $s(G-S) =0$ and $G-S$ is a balanced complete bipartite graph. 
 Let $U$ and $V$ be the two bipartitions of $G-S$. As $G-S$
 is 4.5-tough with respect to $S$ and $c(G-S-U)=c(G-S-V)=|U|=|V|$, there exist $u_1\in U$ and $v_1\in V$
 such that $u_1\sim x_1$ and $v_1\sim y_1$ for distinct $x_1,y_1\in S$.  Let $R_1$
 be a Hamiltonian $(u_1,v_1)$-path of $G-S$ 
  and $Q_1=x_1u_1R_1v_1y_1$. Then $Q_1$ is an $S$-matched basic path-cover with  one component. 
\qed

\begin{LEM}\label{lem:path-partner-system2}
	Let $G$ be a   $(P_4 \cup P_1)$-free graph, and let $S\subseteq V(G)$ be a minimal cutset  for which  $s(G-S) \ge 1$.  
	Suppose that $G$ is 4-tough with respect to $S$. 
	Then $G-S$ has an $S$-matched  basic path-cover $\mathcal{Q}$  such that  each component  $D$ of $G-S$ is covered
	by at most  $\min\{\max\{1, s(D)\},2\}$ components of   $\mathcal{Q}$. 
\end{LEM}

\pf  By Lemma~\ref{lem:path-partner-system1}, $G-S$ has an $S$-matched 
basic path-cover   such that each component  $D$ of $G-S$  is covered 
by $\max\{1, s(D)\}$ components of  the path-cover. 
We choose an  $S$-matched basic 
path-cover  $\mathcal{Q}$   of $G-S$  such that $c(\mathcal{Q})$ is minimized. 

If each component of $G-S$ is covered  by at most two components of $\mathcal{Q}$, 
then we are done. Thus, we 
suppose that some  component $D$ of $G-S$ is covered by  $k$ components 
$Q_1, Q_2, \ldots, Q_k$ of $\mathcal{Q}$, where $k\ge 3$.  This implies that $s(D)\ge 3$. 
Let $S_0\subseteq V(D)$ be 
a maximal scattering set of $D$.  
We suppose $Q_i=x_iu_iR_iv_iy_i$ for each $i\in [1,k]$, where  $R_i:=u_iQ_iv_i$, and $x_i, y_i\in S$.

For distinct $i,j\in [1,k]$, if $E_G(\{x_i, y_i\}, \{ u_j,v_j\}) \ne \emptyset$ or $E_G(\{x_j, y_j\}, \{ u_i,v_i\}) \ne \emptyset$, say $y_i\sim u_j$, then $x_iQ_iy_i u_jQ_jy_j$  and the rest components of $\mathcal{Q}$ form an  $S$-matched
basic path-cover of $G-S$ with fewer components, a contradiction to the choice of $\mathcal{Q}$. 
Thus we assume that there exist  distinct $i,j\in [1,k]$ such that $E_G(\{x_i, y_i\}, \{ u_j,v_j\}) =E_G(\{x_j, y_j\}, \{ u_i,v_i\}) =\emptyset$. 
This particularly  implies that $y_i\sim v_i$ and $v_j \not\sim y_i, v_i$, and $x_j\sim u_j$ and $u_i \not\sim x_j, u_j$. 
As $G$ is $(P_4\cup P_1)$-free and $S$ is a minimal cutset of $G$,  Lemma~\ref{lem:P4-U-P1-cut} implies that 
both $y_i$ and $x_j$ are complete in $G$ to 
all components of $G-S$ other than $D$.  
Thus $y_i$ and $x_j$ have a common neighbor $z$  in $G$ from a component of $G-S$
that is not $D$. Then $v_iy_ix_ju_j$ is an induced $P_4$ in $G$ if $x_j\sim y_i$ and  $v_iy_izx_ju_j$ is an induced $P_5$ in $G$ otherwise. As $G$ is $(P_4\cup P_1)$-free, 
vertices from all  components of $D-S_0$   not containing $v_i$ or $u_j$ are adjacent in $G$
to $y_i$ or $x_j$.   
Let $h\in [1,k]\setminus \{i,j\}$. Then as $\mathcal{Q}$ is an $S$-matched basic path-cover of $G-S$, 
it follows that the vertices  $u_h, v_h$ from  $Q_h$ (recall that $Q_h=x_hu_hR_hv_hy_h$) are from
a component of $D-S_0$  different than the ones  containing vertices $u_i, v_i, u_j, v_j$.  Thus 
$u_h$ and  $v_h$  are adjacent in $G$
to $y_i$ or $x_j$.    
Assume, without loss of generality, that $y_i\sim u_h$. 
Then  $x_iQ_iy_i u_hR_hv_hy_h$  and the rest components of $\mathcal{Q}$ form an  $S$-matched
basic path-cover of $G-S$ with fewer components, a contradiction to the choice of $\mathcal{Q}$. 
\qed

We need the following result by 
H\"aggkvist and  Thomassen from 1982 in the proof of our next lemma. 

\begin{THM}[{\cite[Theorem 1]{Cycles-through-edges}}]\label{thm:cycle-through-ind-edges}
	Let $G$ be a graph and $L$ be a set of $k$ independent edges of $G$, where $k\ge 0$ is an integer. 
	If any two endvertices of edges of $L$ are connected by $k+1$ internally disjoint paths, 
	then $G$ has a cycle containing all edges of $L$. 
\end{THM}

\begin{LEM}\label{lem:cycle-covering-G-S}
	Let $G$ be a $4.5$-tough $(P_4 \cup P_1)$-free graph, and let $S\subseteq V(G)$ be a minimal cutset of $G$. 
	 Then
	\begin{enumerate}[(1)]
		\item $G-S$ has an $S$-matched   path-cover with a single component; and 
		\item $G$ has a cycle covering all vertices of  $G-S$. 
	\end{enumerate}	
\end{LEM}

\pf  
Let $D_1, \ldots, D_\ell$ be all the components of $G-S$, where $\ell \ge 2$ 
is an integer.

When  $\ell \le 3$, 
for $i \in [1,\ell]$, if $s(D_i)\ge 0$ and $D_i$ is not a balanced complete bipartite graph, 
we let $S_i\subseteq V(D_i)$ be a maximal scattering set of $D_i$, and let $z_i$ be 
a  minimal element of $S_i$.
We let $Z$ be the set of all those chosen vertices $z_i$, and let 
$G^*=G-Z$.  

When $\ell \ge 4$, we simply  let $G^*=G$.

We first show that 
$G^*$ is 4-tough with respect to $S$.  Suppose to the contrary that $G^*$ has a cutset $W$
such that $V(D_i-Z)\setminus W\ne \emptyset$ for each $i\in [1,\ell]$ and $\frac{|W|}{c(G^*-W)} <4$. 
For each $i\in [1,\ell]$, if $c(D_i-Z-W) \ge 2$ and $z_i$ exists, we add $z_i$ to $W$; 
if $c(D_i-Z-W) =1$ and $z_i$ exists, then we add  $z_i$  to $D_i-Z$. 
Let $W^*$ be the resulting set of $W$ after adding all the qualified $z_i$'s. Then we have $c(G-W^*)=c(G^*-W)$. 
On the other hand, we have  $|W^*| \le |W|+h$, where $h:= |\{i\in [1,\ell]: c(D_i-Z-W) \ge 2\}|$. 
However, we get $\frac{|W^*|}{c(G-W^*)} \le \frac{|W|+h}{c(G^*-W)} < 4 +\frac{1}{2}=4.5$ (note that $c(G^*-W) \ge 2h$), 
a contradiction to the toughness of $G$. Thus $G^*$ is 4-tough  with respect to $S$.

By Lemma~\ref{lem:path-partner-system2}, $G-S$ has  an $S$-matched basic path-cover $\mathcal{Q}$ such that each  component   $D$   of $G^*-S$ is covered
by at most $\min\{\max\{1, s(D)\},2\}$  components of  $\mathcal{Q}$.   
As $S$ is a cutset of $G$, we know that $c(\mathcal{Q}) \ge 2$. 
Let  $k=c(\mathcal{Q})$ and $Q_1, \ldots, Q_k$ be all the components of $\mathcal{Q}$.
Furthermore, 
we assume that $Q_i=x_iu_iR_iv_iy_i$,  where $x_i, y_i\in S$, and  
$R_i:=u_iQ_iv_i$.  We choose $\mathcal{Q}$ such that the number of components of $\mathcal{Q}$ 
that cover  vertices of  $D_i-Z$  for $i\in [1,\ell]$ is minimized. 
Thus if  there exist distinct $Q_i$ and $Q_j$ that   both contain vertices of $D_i-Z$, then we must have 
\begin{equation}\label{eqn:non-adj-in-H}
E_G(\{x_i, y_i\}, \{u_j,v_i\}) =E_G(\{x_j, y_j\}, \{u_i,v_i\}) =\emptyset. 
\end{equation}

If $\ell \in [2,3]$ and for some  $i\in [1,\ell]$, $D_i-Z$ is covered 
by at least two components of $\mathcal{Q}$, then we must have $s(D_i-Z) \ge 2$
by Lemma~\ref{lem:path-partner-system2}.  Thus $s(D_i) \ge 1$
and so $z_i$ exists.  In this case, we argue that $D_i-Z$ is covered 
by exactly two components of $\mathcal{Q}$, and we use $z_i$
to link the two components into one.   We prove the following claim.

\begin{CLA}\label{claim:H-min-degree0}
	The following statements hold. 
	\begin{enumerate}[(a)]
		\item  For each $i\in [1,\ell]$, $D_i-Z$ is covered 
		by at most two components of $\mathcal{Q}$. 
		If $D_i-Z$ is covered 
		by two  distinct components $Q_h, Q_p$ of $\mathcal{Q}$, 
		then each vertex $z\in \{x_p, y_p, x_h, y_h\}$ is   complete in 
		$G$ to each  $D_j$ with $j\in [1,\ell]$ and $j\ne i$. 
		\item   For each $i\in [1,\ell]$,  if $D_i-Z$ is covered 
		by  exactly one component  $Q_p$  of $\mathcal{Q}$, then 
		there is at most one other component  $Q_h$  of $\mathcal{Q}$
		such that $x_p\not\sim u_h$ and $x_h\not\sim u_p$. 
	\end{enumerate}
\end{CLA}

\pf  For (a),  suppose to the contrary that  for some $i\in [1,k]$,  $D_i-Z$ 
is covered by at least three distinct components 
$Q_h, Q_p, Q_r$  of $\mathcal{Q}$. 
By our assumption, we have $E_G(\{x_p, y_p\}, \{u_h,v_h\}) =E_G(\{x_h, y_h\}, \{u_p,v_p\}) =\emptyset$; and the same conclusion holds with respect to $Q_h, Q_r$ and $Q_p, Q_r$.  
Then $x_p\not\sim u_h$. 
We also have $u_h\not\sim u_p$
by $\mathcal{Q}$ being basic. 
Then by  Lemma~\ref{lem:P4-U-P1-cut} and the fact that $S$ is a minimal cutset of $G$, 
$x_p$ is  complete in $G$ to all components $D_j$ with $j\in [1,\ell]\setminus \{i\}$. 
The same conclusion holds for  $x_h$. 
Thus $x_p$ and $x_h$ have a common neighbor $z$  in $G$ from a component of $G-S$
that is not $D_i$. Then $u_p x_px_hu_h$ is an induced $P_4$ in $G$ if $x_p\sim x_h$  in $G$ and  $u_p x_pzx_hu_h$ is an induced $P_5$ in $G$ otherwise. As $G$ is $(P_4\cup P_1)$-free,  and $u_r\not\sim u_p, u_h$ 
by $\mathcal{Q}$ being basic and $u_r\not\sim z$ as they are from different components of $G-S$,  it follows that $u_r$  is adjacent in $G$ to $x_p$ or $x_h$. 
However, this gives a contradiction to the assumption given in~\eqref{eqn:non-adj-in-H}.  

Thus if  some    $D_i-Z$  
is covered by  exactly two components  $Q_h$ and $Q_p$ of $\mathcal{Q}$,  
then $x_p$ is  complete in $G$ to all components $D_j$ with $j\in [1,\ell]\setminus \{i\}$. 
Similarly, each of $y_p, x_h, y_h$ 
is  complete in $G$ to all components $D_j$ with $j\in [1,\ell]\setminus \{i\}$.

For (b),  suppose that $x_p\not\sim u_h$ and $x_h\not\sim u_p$. 
Then we also have $u_h\not\sim u_p$ in $G$
by $\mathcal{Q}$ being   basic. 
 As $S$ is a minimal cutset of $G$, by  Lemma~\ref{lem:P4-U-P1-cut}, 
$x_p$ is  complete in $G$ to all other components of $G-S$ other than the one containing $u_h$.  The corresponding conclusion holds also for $x_h$. 
By Claim~\ref{claim:H-min-degree0}(a), it must be the case that 
some $D_j-Z$  for $j\in [1,\ell]$ with $j\ne i$ that contains 
  $u_h$  is covered by a single component of $\mathcal{Q}$.  
 Thus, there is at most one other component  $Q_h$  of $\mathcal{Q}$
 such that $x_p\not\sim u_h$ and $x_h\not\sim u_p$.    
\qed 

Now for each $i\in [1,\ell]$, if $D_i-Z$ is covered by exactly two components 
$Q_h$ and $Q_p$ of $\mathcal{Q}$, then we know that $z_i$ exists. 
Then we link the two components $Q_h$ and $Q_p$ of $\mathcal{Q}$ into one using 
$z_i$. Precisely, we  let  $x_hu_hQ_hv_hz_iu_pQ_pv_py_p$ be 
 a new basic path-cover of $D_i$.  Note that by Claim~\ref{claim:H-min-degree0}, 
 $x_h$ and $y_p$ are complete in $G$ to all other components  
 of $G-S$. For notational simplicity, we continue to use $\mathcal{Q}$ to denote the new $S$-matched basic path cover, and we retain the same notation as before. The only points to keep in mind are the following.

\begin{CLA}\label{claim:Q-property}
	\begin{enumerate}[(a)]
		\item If $z_i$ exists for $D_i$ and $z_i$ is not a vertex of $\mathcal{Q}$ for some $i\in [1,\ell]$, then $D_i-z_i$ is covered by one single 
		component of $\mathcal{Q}$.   
		\item If $z_i$ exists for $D_i$ and $z_i$ is  a vertex of $\mathcal{Q}$  for some  $i\in [1,\ell]$, then   the $S$-endvertices 
		of the component of $\mathcal{Q}$ covering $D_i$ are complete in $G$
		to all other $D_j$ with $j\in [1,\ell]$ and $j\ne i$. 
	\end{enumerate}
\end{CLA}

We now construct an auxiliary graph  $H$  and use that to demonstrate the existence of a single 
path or cycle that covers all vertices of $G-S$.  The graph $H$ is constructed as follows. 
Its  vertices are  $x_1, y_1, \ldots, x_k, y_k$, and $E(H)$ consists of 
$x_1y_1, \ldots, x_ky_k$, and additionally a vertex $x$ is adjacent in $H$ to a vertex $y$
if $x$ is adjacent in $G$ to the partner of $y$ in $\mathcal{Q}$ or $y$ is adjacent in $G$ to the partner of $x$ in $\mathcal{Q}$. 
By this construction,  $H$ is a graph on $2k$ vertices.      By Claim~\ref{claim:H-min-degree0},  $H$ satisfies the following property.

\begin{CLA}\label{claim:H-min-degree}
We have $\delta(H) \ge 2k-3$. Furthermore, if $d_H(z)=2k-3$ for some $z\in \{x_1,y_1, \ldots, x_k,y_k\}$, then 
the two non-neighbors of $x$ are $x_j,y_j$ for some $j\in [1,k]$. 
\end{CLA}

 When $k\ge 5$,  we show that $H$ is $(k+1)$-connected. For otherwise, $H$ has a cutset $W$ of size at most $k$. 
 As each vertex of $H$ has degree at least $2k-3$ in $H$, it follows that each component of $H-W$ contains at most two 
 vertices.  
On the other hand, by  $\delta(H) \ge 2k-3$,  we know that each component of $H-W$ has at least $2k-2-|W|$ vertices. 
Thus  $2\ge 2k-2-|W|$, giving $|W| \ge 2k-4$. This combined with  $|W|  \le k$, gives $k\le 4$, a contradiction.  
Thus $H$ is $(k+1)$-connected.  By Theorem~\ref{thm:cycle-through-ind-edges}, $H$ contains a cycle  $C$ and so also 
a path $P$ such that $C$ and $P$ contains all the edges  $x_1y_1, \ldots, x_ky_k$.  
For each $i\in [1,k]$, we replace $x_iy_i$ on $C$ and $P$ by $Q_i$. 
For an edge $xy\in E(C) \cup E(P)$ such that $x$ and $y$ are from different components of $\mathcal{Q}$, 
we let $x'$ and $y'$ be respectively the partners of $x$ and $y$ in $\mathcal{Q}$.  
By the construction of $H$, we know that $xy'\in E(G)$ or $yx'\in E(G)$. 
We then  replace $xy$ by one edge in $\{xy', yx'\}\cap E(G)$. After these replacements, the resulting cycle of  $C$ 
is a cycle  covering all vertices of $G^*-S$, and the resulting path of 
  $P$ is an $S$-matched basic path-cover of $G^*-S$ with one single component. 
For any $z\in Z$ such that $z$ is not a vertex of  $\mathcal{Q}$,   we  let $D$
be the component of $G-S$ containing $z$. Then we know that $s(D) \ge 0$
and $D$ is not a balanced complete bipartite graph by the existence of $z$. 
By Claim~\ref{claim:Q-property}(a) and the construction of $P$ and $C$, 
we know that the longest segment of each $P$ and $C$ starts 
and ends with an $S$-vertex containing  a vertex of $D$ is an $S$-matched 
basic path-cover with one single component.  By Lemma~\ref{lem:inserting},
we can insert $z$ in the middle of that segment.  After this modification 
for all such vertices $z$, we have constructed a desired path and cycle 
covering vertices of $G-S$.

When $k=4$,   if $H$ is $(k+1)$-connected, then we can construct a desired cycle or path covering vertices
 of $G-S$ the same way as above.  Thus we assume that $H$ is not $(k+1)$-connected. Then by $\delta(H) \ge 2k-3$, 
 it follows that $H$ has a cutset  $W$ consisting 
of exactly 4 vertices for which $H-W$ has exactly two components that each consists of an edge of the form 
$x_iy_i$ for some $i\in [1,k]$.  (For a vertex $x$ of $H$ that has two non-neighbors from $V(H)\setminus \{x\}$, the two non-neighbors 
form an edge from $\{x_1y_1, \ldots, x_ky_k\}$ by Claim~\ref{claim:H-min-degree}).  Furthermore,  the subgraph of $H$ 
induced by the edges between $W$ and $V(H)\setminus W$ is a complete bipartite graph by $\delta(H) \ge 5$. 
Assume, without loss of generality that $x_1, y_1, x_2, y_2 \in W$ and $x_3y_3$
and $x_4y_4$ are respectively the two components of $H-W$. Then   $x_1y_1x_3y_3x_2y_2x_4y_4$ and $x_1y_1x_3y_3x_2y_2x_4y_4x_1$
are respectively a path and a cycle containing $x_1y_1, \ldots, y_4y_4$ in $H$. 
Then we can construct a desired cycle and  path covering vertices
of $G-S$ the same way as the case $k\ge 5$.

  Thus we are only left to construct a desired path and cycle when $k\in [2,3]$.    
  If  vertices of $D_i-Z$ is covered by two components of $\mathcal{Q}$, then we must have $s(D_i-Z) \ge 2$ by Lemma~\ref{lem:path-partner-system2}. This implies that the vertex $z_i$ 
  exists.  Thus by Claim~\ref{claim:Q-property}(a), we have $k=\ell$. 
  By renaming  the components $D_i$   if necessary, we assume that $Q_i$
  covers all vertices of $D_i-Z$ for each $i\in [1,\ell]$. 
  We consider two cases.  
  
   {\bf \noindent  Case 1: $k =\ell =3$.}
   
    If $\delta(H) \ge 4$, then $H$ is 4-connected and so we can find a desired path and cycle 
   the same way as in the $k\ge 5$ case. Thus we assume that $\delta(H) =3$ 
   by Claim~\ref{claim:H-min-degree}. Furthermore, without loss of generality, 
   we assume that $x_1 \not\sim x_2,y_2$ in $H$. This implies that  $E_G(x_1, \{u_2,v_2\}) =E_G(u_1,\{x_2,y_2\}) =\emptyset$.  Application of  Lemma~\ref{lem:P4-U-P1-cut} implies that 
   each of $x_1, x_2, y_2$ is  complete in $G$ to all vertices of $D_3$. 
   If $y_1$ is adjacent in $H$ to one of $x_2,y_2$, then  $H$ has a cycle containing $x_1y_1, x_2y_2$ and $x_3y_3$, and so we can find a desired path and cycle 
   the same way as in  the $k\ge 5$ case also. 
   Thus $y_1$ is  complete in $G$ to all vertices of $D_3$ as well.

  As $S$ is a minimal cutset of $G$,  $y_1$ has in $G$ a neighbor $w_2$ from $V(D_2)$.  
  We construct a dersied path and cycle in each of the following cases.

  If $w_2\in \{u_2, v_2\}$, say $w_2=u_2$,  then we let 
  $P^*=x_1Q_1v_1y_1u_2Q_2y_2u_3Q_3v_3y_3$ and $C^*=x_1Q_1v_1y_1u_2Q_2y_2u_3Q_3v_3x_1$. 
  Thus $w_2\not\in \{u_2, v_2\}$.

  If $s(D_2) \le -1$, then $D_2$ has a Hamiltonian $(w_2,v_2)$-path $R_2^*$. Let 
  $$P^*=x_1Q_1v_1y_1w_2R^*_2v_2y_2u_3Q_3v_3y_3 \quad \text{and} \quad C^*=x_1Q_1v_1y_1w_2R^*_2v_2y_2u_3Q_3v_3x_1.$$
  
  If $s(D_2) =0$ and $D_2$ is a balanced complete bipartite graph, then $u_2$ and $v_2$ are from different bipartitions of $D_2$. 
  Thus there is in $D_2$ a Hamiltonian path $R_2^*$ from $w_2$ to exactly one of $u_2$ and $v_2$, say to $v_2$ without loss of generality. 
  Then we let $$P^*=x_1Q_1v_1y_1w_2R^*_2v_2y_2u_3Q_3v_3y_3 \quad \text{and} \quad C^*=x_1Q_1v_1y_1w_2R^*_2v_2y_2u_3Q_3v_3x_1.$$ 
  
  For  the both cases above, 
  if $z_1$ or $z_3$ exist, then we can respectively insert them within the segments $u_1Q_1v_1$ or  $u_3Q_3v_3$ of both $P^*$ and $C^*$  
  by Lemma~\ref{lem:inserting} 
  to get a desired path and cycle. 
  
  Thus we assume that  $s(D_2) \ge 0$ and $D_2$ is not a balanced complete bipartite graph. Then the vertex $z_2$ exists. Also by the assumption that 
  $E_G(u_1,\{x_2,y_2\}) =\emptyset$ and Claim~\ref{claim:Q-property}(b), we know that $z_2$ is not a vertex of $Q_2$. 
  \begin{itemize}
  	\item If $w_2=z_2$, then as $z_2\sim u_2,v_2$, we let $$P^*=x_1Q_1v_1y_1z_2u_2Q_2y_2u_3Q_3v_3y_3 \quad \text{and} \quad  C^*=x_1Q_1v_1y_1z_2u_2Q_2y_2u_3Q_3v_3x_1.$$ If $z_1$ or $z_3$ exist, then we can respectively insert them within the segments $u_1Q_1v_1$ or  $u_3Q_3v_3$ of both $P^*$ and $C^*$ 
  	to get a desired path and cycle  by Lemma~\ref{lem:inserting}. 
  	\item Thus we assume that $w_2\ne z_2$. Since $w_2\not\in \{u_2, v_2\}$ also, $w_2$ is an internal vertex of $u_2Q_2v_2$. 
  	Let $w_2^-$ and $w_2^+$ be respectively the two neighbors of $w_2$ on $u_2Q_2v_2$, where $w_2^-$ lies on $u_2Q_2w_2$. 
  	If $z_2$ is adjacent in $G$ to one of $w_2^-$ and $w_2^+$, say $w_2^-$, then we 
  	let 
  	\begin{eqnarray*}
  		P^* &=&x_1Q_1v_1y_1w_2Q_2v_2z_2w_2^-Q_2u_2x_2u_3Q_3v_3y_3, \\
  		C^*&=&x_1Q_1v_1y_1w_2Q_2v_2z_2w_2^-Q_2u_2x_2u_3Q_3v_3x_1.
  	\end{eqnarray*}
  	If $z_1$ or $z_3$ exist, then we can respectively  insert them within the segments $u_1Q_1v_1$ or $u_3Q_3v_3$ of both $P^*$ and $C^*$ 
  	to get a desired path and cycle.
  	\item Thus we assume that $z_2\not\sim w_2^-, w_2^+$. This implies that both $w_2^-$ and $w_2^+$ are minimal elements of $S_2$ in $D_2$. 
  	Then we let 
  	\begin{eqnarray*}
  		P^*&=&x_1Q_1v_1y_1w_2Q_2v_2w_2^-Q_2u_2x_2u_3Q_3v_3y_3, \\ 
  		C^*&=&x_1Q_1v_1y_1w_2Q_2v_2w_2^-Q_2u_2x_2u_3Q_3v_3x_1. 
  	\end{eqnarray*}
  	Now, if exist,  we insert  $z_1$, $z_2$ or $z_3$ respectively within  segments  $u_1Q_1v_1$, $w_2Q_2v_2w_2^-Q_2u_2$, or  $u_3Q_3v_3$
  	of $P^*$ and $C^*$ to get the desired path and cycle. 
  \end{itemize}

     {\bf \noindent  Case 2: $k =\ell=2$.}

 By Claim~\ref{claim:Q-property}(b), if $z_1$ or $z_2$ exists, 
 we may assume that none of them is a vertex of $\mathcal{Q}$ (as otherwise, we can easily construct a desired path and cycle). 
  We first make the following claim. 
   \begin{CLA}\label{claim:assumptions}
   We can make the 
   following assumptions: 
    \begin{enumerate}[(1)]
   	\item $y_1$ has in $G$ a   neighbor $w_2$  from $V(D_2)\setminus \{v_2\}$. Furthermore, if $D_2$ is  a balanced complete bipartite graph, then 
   	$w_2$ and $v_2$ are from different bipartitions of $D_2$; 
   	
   	\item $y_2$ has in $G$ a neighbor  $w_1$ from $V(D_1)\setminus \{v_1\}$. Furthermore,  if $D_1$ is  a balanced complete bipartite graph, then 
   	$w_1$ and $v_1$ are from different bipartitions of $D_1$. 
   \end{enumerate}
   \end{CLA}

    \proof[Proof of Claim~\ref{claim:assumptions}]
   We suppose to the contrary, and without loss of generality,  that $w_2=v_2$ when $D_2$ is not a balanced complete bipartite graph, and 
   $w_2$ and $v_2$ are from  the same bipartition of $D_2$ when $D_2$ is a balanced complete bipartite graph. 
   
   If $x_2$ has in $G$ a neighbor from $V(D_1)$ that is not $v_1$ when $D_1$
   is not a balanced complete bipartite graph, and  is not in the same bipartition as $v_1$ when $D_1$ is a balanced complete bipartite graph, then 
   we can just exchange  the labels of $u_2$ and $v_2$ and that of $x_2$ and $y_2$ 
   in getting our desired assumption. 
   
   Thus we assume that  $x_2$ has in $G$ a neighbor from $V(D_1)$, and the neighbor   is  only  $v_1$ when $D_1$
   is not a balanced complete bipartite graph, and  is   in the same bipartition as $v_1$ when $D_1$ is a balanced complete bipartite graph. 
   We then consider a neighbor  $w$ of $x_1$ in $G$ from $V(D_2)$. 
   If $w=u_2$, then  let  $P^*=x_1u_1Q_1v_1y_1v_2Q_2u_2x_2$  and $C^*=x_1u_1Q_1v_1y_1v_2Q_2u_2x_1$.  
   If $z_1$ or $z_2$ exist, by Lemma~\ref{lem:inserting}, we can insert them respectively  in the segments $u_1Q_1v_1$ or  $v_2Q_2u_2$ 
   of $P^*$ and $C^*$ and get our desired path and cycle. 
   Thus 
   we assume that $w\ne u_2$. 
   If $D_2$ is a balanced complete bipartite graph and $w$ and $u_2$ are from the same bipartition of $D_2$, 
   then  $w$ and $v_2$ are from different bipartitions of $D_2$. We let $R_2^*$ be a  Hamiltonian $(w,v_2)$-path of $D_2$, and let $Q_2^*=wR_2^*v_2y_2$. 
   Let $P^*=y_1v_1Q_1u_1x_1wR_2^*v_2y_2$  and $C^*=x_1u_1Q_1v_1y_1v_2Q^*_2wx_1$. 
   If $z_1$ or $z_2$ exist, we can insert them respectively  in the segments $u_1Q_1v_1$ or $wR_2^*v_2$
   of $P^*$ and $C^*$ and get our desired path and cycle. 
   Thus  we assume that $w\ne u_2$, and when $D_2$ is a balanced complete bipartite graph then  $w$ and $u_2$ are from different  bipartitions of $D_2$. 
   Then 
   exchanging   the labels of  $u_1$ and $v_1$, of $x_1$ and $y_1$, of $u_2$ and $v_2$,    and of $x_2$ and $y_2$ 
   gives our desired assumption. 
   \qed 
   
   If  $s(D_1) \le -1$ or $s(D_1)=0$ and $D_1$ 
   is a balanced complete bipartite graph (so the vertex $z_1$ does not exist), then we let $R_1^*$ be a Hamiltonian $(w_1,v_1)$-path of $D_1$. 
   If  $s(D_2) \le -1$ or $s(D_2)=0$ and $D_2$ is a balanced complete bipartite graph (so the vertex $z_2$ does not exist), 
   then we let $R_2^*$ be a Hamiltonian $(w_2,v_2)$-path of $D_2$.   
   We now construct a desired path and cycle according to the size of $Z$. 
   
   If $Z=\emptyset$, then the above two cases happen and we let 
   \begin{eqnarray*}
   	P&=&x_1u_1Q_1v_1y_1w_2R_2^*v_2y_2, \\
   	C&=& w_1R_1^*v_1y_1w_2R_2^*v_2y_2w_1, 
   \end{eqnarray*}
   which are respectively our desired path and cycle. 
   
   Next we consider $|Z|=1$, and by symmetry, we assume that $Z=\{z_1\}$. If $w_1=u_1$, then we can construct $P$ and $C$ the same 
   as above, but insert $z_1$ in the segment $u_1Q_1v_1$ of $P$ and $C$ to get our desired path and cycle. Thus we assume that $w_1\ne u_1$. 
   Let $P^*=x_1u_1Q_1v_1y_1w_2R_2^*v_2y_2$. Then a desired path is obtained 
   from $P^*$ by inserting $z_1$ in the segment $u_1Q_1v_1$ of $P^*$.  Now we  construct a desired cycle in this case. 
   As $w_1\ne v_1$ by our assumption, $w_1$ is an internal vertex of $Q_1$. 
   Let $w_1^-$ and $w_1^+$ be respectively the two neighbors of $w_1$ on $Q_1$, where $w_1^-$ lies on $u_1Q_1w_1$. 
   If $z_1 \sim  w_1^+$, then  $C:=w_1Q_1u_1z_1w_1^+Q_1v_1y_1w_2R_2^*v_2y_2w_1$ is a desired cycle.  
   If $z_1 \not\sim  w_1^+$, then $w_1^+$ is also a minimal element of $S_1$.  Let $C^*=w_1Q_1u_1w_1^+Q_1v_1y_1w_2R_2^*v_2y_2w_1$. 
   Then a desired cycle  is obtained 
   from $C^*$ by inserting $z_1$ in the segment $w_1Q_1u_1w_1^+Q_1v_1$ of $C^*$.
   
   Lastly, we assume that $Z=\{z_1,z_2\}$ and consider three subcases as follows.

   If $w_1=z_1$ and $w_2=z_2$, then we let 
   \begin{eqnarray*}
   	P^*&=&x_1u_1Q_1v_1y_1w_2u_2Q_2v_2y_2, \\
   	C&=& w_1u_1Q_1v_1y_1w_2u_2Q_2v_2y_2w_1. 
   \end{eqnarray*}
   Then $C$ is our desired cycle, and a desired path is obtained from $P^*$ by inserting $z_1$ in the segment $u_1Q_1v_1$ of $P^*$. 
   
   For the second subcase, by symmetry, we assume that  $w_1 \ne z_1$ and $w_2= z_2$.  We let $P^*=x_1u_1Q_1v_1y_1w_2u_2Q_2v_2y_2$. 
   Then we insert $z_1$ in the segment $u_1Q_1v_1$ of $P^*$  
   in getting our desired path. 
   If  $w_1=u_1$,  then we  let  $C^*=u_1Q_1v_1y_1w_2u_2Q_2v_2y_2u_1$ 
   and insert $z_1$ in the segment $u_1Q_1v_1$ of  $C^*$  
   in getting our desired   cycle. 
   Thus  we assume  $w_1\ne u_1$.  As also $w_1\ne v_1$ by Claim~\ref{claim:assumptions}, we know that $w_1$ is an internal vertex of $u_1Q_1v_1$. 
   Let $w_1^-$ and $w_1^+$ be respectively the two neighbors of $w_1$ on $Q_1$, where $w_1^-$ lies on $u_1Q_1w_1$. 
   If $z_1 \sim  w_1^+$, then   $C:=w_1Q_1u_1z_1w_1^+Q_1v_1y_1w_2u_2Q_2v_2y_2w_1$ is our desired cycle.  
   If $z_1 \not\sim  w_1^+$, then $w_1^+$ is also a minimal element of $S_1$. We  let $C^*=w_1Q_1u_1w_1^+Q_1v_1y_1w_2u_2Q_2v_2y_2w_1$ 
   and  insert $z_1$ in the segment $w_1Q_1u_1w_1^+Q_1v_1$ of  $C^*$  
   in getting our desired   cycle.

   Lastly, we consider  $w_1\ne z_1$ and $w_2\ne z_2$. Note that $w_2\ne v_2$ by Claim~\ref{claim:assumptions}. 
   Let $w_2^+$ be the neighbor of $w_2$ lying on the path $w_2Q_2v_2$. 
   If $z_2\sim w_2^+$, then  we let $R_2^*=w_2Q_2u_2z_2w_2^+Q_2v_2y_2$. 
   Thus we assume that  $z_2\not\sim w_2^+$. This implies that $w_2^+$ is also a minimal element of $S_2$ in $D_2$. 
   Then we let $R_2^*$ be obtained from $w_2Q_2u_2w_2^+Q_2v_2y_2$ by inserting $z_2$.  Let 
   $P^* =x_1u_1Q_1v_1y_1w_2R_2^*v_2y_2$. Then we insert $z_1$ in the segment $u_1Q_1v_1$ of $P^*$  
   in getting our desired path.  
   In the same way as above, we can also find a Hamiltonian $(w_1,v_1)$-path $R_1^*$ of $D_1$ (containing the vertex $z_1$). 
   Then $C=w_1R_1^*v_1w_2R_2^*v_2y_2w_1$ is our desired cycle. 
   \qed 

\subsection{Construct a Hamiltonian  cycle when a suitable cutset is given}

Let $\oC$ be an oriented cycle. For $x\in V(C)$, denote the immediate  successor of $x$ by $x^{+}$ and the  immediate predecessor of $x$ by $x^{-}$ following the 
orientation of $C$. For $u,v\in V(C)$, $u\oC v$ denotes the segment  of $C$ starting with $u$, following $C$ in the orientation, and ending at $v$. Likewise, $u \iC v$ is the opposite segment of $C$ with ends $u$ and $v$. We assume all cycles in consideration afterwards are oriented in the clockwise direction.

\begin{LEM}\label{lem:vertex-insetting}
	Let $t>0$ and $G$ be a $t$-tough  $n$-vertex  graph with a non-Hamiltonian cycle  $C$. For a connected subgraph  $H$ of $G-V(C)$, if 
	 $|N_G(H,C)|> \frac{n}{t+1}-1$, then we can extend $C$ to a cycle $C^*$ 
such that  $V(C) \subseteq V(C^*)$ and $V(C^*)\cap V(H) \ne \emptyset$. 
\end{LEM}

\pf Let $v_1, \ldots, v_k$ be all the neighbors of  vertices of $H$ on $C$, and we assume that these vertices appear in
the order $v_1, \ldots, v_k$ along $\oC$, where $k\ge 1$ is an integer.   If $v_iv_{i+1} \in E(C)$ for some $i$, where the indices are taken modulo $k$, 
then we let  $v_i^*, v^*_{i+1} \in V(H)$ such that $v_i^*\sim v_i$ and  $v^*_{i+1} \sim v_{i+1}$, and  let $P$ be a $(v^*_{i}, v^*_{i+1})$-path in $H$. 
Now  $C^*=v_{i+1}\oC v_i v_i^*Pv^*_{i+1}v_{i+1}$ is a desired cycle. Thus we assume that no two vertices among $v_1, \ldots, v_k$
are consecutive on $C$. If for some $i,j\in [1,k]$, say without loss of generality, that $i<j$, we have $v_i^+\sim v_j^+$, 
then we let  $v_i^*, v^*_{j} \in V(H)$ such that $v_i^*\sim v_i$ and  $v^*_{j} \sim v_{j}$, and  let $P$ be a $(v^*_{i}, v^*_{j})$-path in $H$. 
Now $C^*=v_j^+\oC v_i v_i^*Pv_j^*v_{j}\iC v_i^+v_j^+$ is a desired cycle.  Thus we assume that 
$\{v_1^+, \ldots, v_k^+\}$ is an independent set of $G$, and $x\not\sim v^+_i$ for any $i\in [1,k]$ and any $x\in V(H)$. 
Let $x\in V(H)$. 
Then  $W:=\{x, v^+_1, \ldots, v^+_k\}$ is an independent set in $G$. 
However, $2 \le |W|=k+1=d_G(x,C)+1>\frac{n}{t+1}$ and so $\frac{|V(G)\setminus W|}{|W|}<t$, 
a contradiction to $G$ being $t$-tough. 
\qed 

\begin{LEM}\label{lem:S-vertex-large-degree-in-G-S}
Let   $G$ be a $4.5$-tough $(P_4\cup P_1)$-free $n$-vertex graph, and $S\subseteq V(G)$  be a cutset of $G$. 
For any subset $S_0\subseteq S$, 
if there is an ordering  ``$<$'' of vertices of $S_0$: $x_1< x_2< \ldots < x_{s_0}$, where $s_0:=|S_0|$, 
such that $d_G(x_i, (V(G)\setminus S)\cup \{x_{1}, \ldots, x_{i-1}\})
> \frac{n}{t+1}-1$, then $G$ has a cycle containing all vertices of $(V(G)\setminus S)\cup S_0$. 
\end{LEM}

\pf  By removing  vertices of $S$  to $G-S$ if necessary, we assume that $S$ is a minimal cutset of $G$. 
Note that removal of vertices  preserves the degree condition for the remaining   vertices of $S_0$.  
Applying Lemma~\ref{lem:cycle-covering-G-S}, we let $C$ be a cycle of $G$ 
that covers all the vertices of $G-S$. Let $S_1=S_0\setminus V(C)$. If $S_1=\emptyset$,
then $C$ is a desired cycle already. Thus we assume that  $S_1 \ne \emptyset$. 
Let $s_1=|S_1|$ and   $S_1=\{y_1, \ldots, y_{s_1}\}$.  We further assume 
that the labels of the vertices of $S_1$ are chosen so that  $y_1<y_2<\ldots < y_{s_1}$. 
Applying Lemma~\ref{lem:vertex-insetting} with $H=y_1$, we find a cycle $C_1$
such that $V(C_1)=V(C) \cup \{y_1\}$. Now for each $i\in [2,s_1]$, 
we apply Lemma~\ref{lem:vertex-insetting} with $H=y_i$ and cycle $C_{i-1}$,
we get a cycle $C_i$ such that $V(C_i)=V(C_{i-1})\cup \{y_i\}$. Then $C_{s_1}$
is our desired cycle. 
\qed

\begin{THM}\label{thm:two-large-component-case}
	Let $G$ be a $4.5$-tough $(P_4\cup P_1)$-free graph on $n\ge 3$ vertices, and let $S$ be a cutset of $G$.   If  $G-S$ 
	has  one  component of order at least $\frac{2n}{t+1}$ and the total order of the  others  is at least $\frac{2n}{t+1}$,  then $G$ is Hamiltonian. 
\end{THM}

\pf Let $D_1, \ldots, D_\ell$ be all the components of $G-S$, where $\ell \ge 2$ 
is an integer.   Without loss of generality, we assume that  $|V(D_1)| \ge \frac{2n}{t+1}$. 
If there is $x\in S$ such that  $N_G(x, D_1)=\emptyset$, then we move $x$ out from $S$. 
Also, if $x\in S$ is   connected in $G$ to  none of the components $D_2, \ldots, D_\ell$, we also 
move $x$ out of $S$. 
Note that $G-(S\setminus \{x\})$ still has    one  component of order at least $\frac{2n}{t+1}$ and the others  of total order at least $\frac{2n}{t+1}$. 
Thus we assume that every vertex of $S$ has in $G$ a neighbor from $D_1$, and is connected to 
at least two components of $G-S$.

We consider two cases regarding whether or not $c(G-S) \ge 3$. 

{\bf \noindent Case 1: $c(G-S)\ge 3$.}

\begin{CLA}\label{claim:S-vertex-adjacency}
	Let $x\in S$.  If $V(D_1)\not\subseteq N_G(x)$, then $x$ is complete to each component $D_i$ with $i\in [2,\ell]$. 
	As a consequence, we have $d_{G}(x, G-S) \ge \frac{2n}{t+1}$ for each $x\in S$. 
\end{CLA}

\proof[Proof of Claim~\ref{claim:S-vertex-adjacency}]
Let $u\in N_{G-S}(x)\setminus V(D_1)$. Assume, without loss of generality, that $u\in V(D_2)$.  As $D_1$ is connected, there is an edge in $D_1$ between  $ N_G(x, D_1)$
and $V(D_1)\setminus N_G(x)$. Thus we 
can  choose $vw\in E(D_1)$ such that $xv\in E(G)$ but $xw\not\in E(G)$. Then $uxvw$ is an induced $P_4$ in $G$. 
As $G$ is $(P_4\cup P_1)$-free, we must have $\bigcup_{i=3}^s V(D_i)\subseteq N_G(x)$. 
Now with $D_3$ in the place of $D_2$, by the same argument as above, we conclude that $V(D_2)\subseteq N_G(x)$. 
Therefore $x$ is complete to  each component $D_i$ with $i\in [2,\ell]$.  The consequence part 
of the statement is clear by the assumption that  $\sum_{i=2}^\ell |V(D_i)| \ge \frac{2n}{t+1}$. 
\qed

Now by Claim~\ref{claim:S-vertex-adjacency} and  Lemma~\ref{lem:S-vertex-large-degree-in-G-S}, $G$ has a Hamiltonian cycle.

{\bf \noindent Case 2: $c(G-S)=2$.}  

By moving a vertex of $S$ to $D_1$ or $D_2$ if necessary, we may assume that  $S$ is a minimal cutset of $G$. 
By the assumption of this theorem, we have $|V(D_i)| \ge \frac{2n}{t+1}$ for each $i\in [1,2]$. 
 Let $S_0=\{x\in S: |N_G(x)\cap V(D_1 \cup D_2)|< \frac{n}{t+1}\}$.  By the definition of $S_0$, 
 for every $x\in S_0$, we have $V(D_i) \setminus N_G(x) \ne \emptyset$ for each $i\in [1,2]$.

\begin{CLA}\label{claim:S0-vertex-adjacency}
For any distinct $x,y\in S_0$,  we have  $N_{G}(x, D_1) \setminus  N_G(y, D_1) =\emptyset$ or $N_{G}(y, D_1) \setminus  N_G(x, D_1) = \emptyset$.  
\end{CLA}
\proof[Proof of Claim~\ref{claim:S0-vertex-adjacency}]  
As   $V(D_i) \setminus N_G(x) \ne \emptyset$ for each $i\in [1,2]$, 
we let $u,v\in V(D_1)$ such that $uv\in E(D_1)$, $x\sim u$, and $x\not\sim v$, and let $w\in N_G(x, D_2)$. 
Then $wxuv$ is an induced $P_4$ in $G$. As $G$ is $(P_4\cup P_1)$-free, we know that $w$ is adjacent in $G$
to every vertex of  $V(D_2) \setminus N_G(x)$. Similarly, by exchanging the roles of $D_1$ and $D_2$ and repeating the same argument, we know that   every neighbor of $x$ in $D_1$ is adjacent in 
$G$ to every vertex of  $V(D_1) \setminus N_G(x)$. The same assertions hold for $y$. 

Assume first that  $x\not\sim y$. 
If  $N_{G}(x, D_2) \setminus  N_G(y, D_2) \ne \emptyset$ and $N_{G}(y, D_2) \setminus  N_G(x, D_2) \ne \emptyset$, 
we choose $u\in N_{G}(x, D_2) \setminus  N_G(y, D_2)$ and $v\in N_{G}(y, D_2) \setminus  N_G(x, D_2) $. 
By the argument in the first paragraph of this proof, we have $uv\in E(D_2)$. 
Then $xuvy$ is an induced $P_4$ in $G$. As $G$ is $(P_4\cup P_1)$-free, 
we know that every vertex of $V(D_1)$ is adjacent in $G$ to $x$ or $y$, and so $\max\{d_G(x, D_1), d_G(y, D_1)\} \ge  \frac{1}{2}|V(D_1)| \ge \frac{n}{t+1}$, a contradiction to $x,y\in S_0$. Thus we must have $N_{G}(x, D_2) \setminus  N_G(y, D_2) =\emptyset$ or $N_{G}(y, D_2) \setminus  N_G(x, D_2) = \emptyset$. Assume, without loss of generality, that $N_{G}(y, D_2) \setminus  N_G(x, D_2) = \emptyset$. 
Thus $N_{G}(y, D_2)\subseteq  N_G(x, D_2)$. In particular, this implies that every vertex of  $V(D_2)\setminus  N_G(x, D_2)$ is  in $G$ 
a common non-neighbor of  $x$ and $y$. 

If $N_{G}(x, D_1) \setminus  N_G(y, D_1) \ne \emptyset$ and $N_{G}(y, D_1) \setminus  N_G(x, D_1) \ne \emptyset$, 
we choose $u\in N_{G}(x, D_1) \setminus  N_G(y, D_1)$ and $v\in N_{G}(y, D_1) \setminus  N_G(x, D_1) $. 
By the argument in the first paragraph of this proof, we have $uv\in E(D_2)$.   Then $xuvy$ is an induced $P_4$ in $G$, which together with a vertex of $V(D_2)\setminus  N_G(x, D_2)$ 
form an induced $P_4\cup P_1$ in $G$,  a contradiction. Thus  we must have $N_{G}(x, D_1) \setminus  N_G(y, D_1) =\emptyset$ or $N_{G}(y, D_1) \setminus  N_G(x, D_1) = \emptyset$.  

Assume then that $x\sim y$. If $N_{G}(x, D_1) \setminus  N_G(y, D_1)  \ne \emptyset$, then we let 
$u\in N_{G}(x, D_1) \setminus  N_G(y, D_1)$ and $v\in V(D_1)\setminus (N_G(x, D_1)\cup N_G(y, D_1))$. 
By the argument in the first paragraph of this proof, we have $uv\in E(D_1)$. 
Then $yxuv$ is an induced $P_4$ in $G$. This implies that every vertex of $D_2$ is 
adjacent in $G$ to $x$ or $y$. Thus  $\max\{d_G(x, D_2), d_G(y, D_2)\} \ge  \frac{1}{2}|V(D_2)| \ge \frac{n}{t+1}$, a contradiction to $x,y\in S_0$.  
Thus $N_{G}(x, D_1) \setminus  N_G(y, D_1)  = \emptyset$. (In fact, in this case, we also have  $N_{G}(y, D_1) \setminus  N_G(x, D_1)  = \emptyset$ 
and so $N_{G}(x, D_1) =  N_G(y, D_1)$.)
\qed

Let $x\in S_0$ such that $d_{G}(x, D_1)$ is largest among that of all vertices of $S_0$. 
Then for any $y\in S_0$ with $y\ne x$, we have  $N_G(y, D_1) \subseteq N_G(x, D_1)$. 
Note that $|N_G(x, D_1)|<\frac{n}{t+1}$ and for any $z\in N_G(x, D_1)$, 
we have $d_G(z, V(D_1)\setminus N_G(x, D_1))>\frac{n}{t+1}$ by the argument in the first paragraph of this proof. 
Now we let $S^*=(S\setminus S_0)\cup  N_G(x, D_1)$. Then $S^*$ is a cutset  of $G$
with the property that every vertex of  $N_G(x, D_1)$ has more than $\frac{n}{t+1}$
neighbors from $V(G)\setminus S^*$, and every vertex of $S^*\setminus N_G(x, D_1)$ 
has at least $\frac{n}{t+1}$ neighbors from $(V(G)\setminus S^*) \cup N_G(x,D_1)$. 
Now by Lemma~\ref{lem:S-vertex-large-degree-in-G-S}, $G$ has a Hamiltonian cycle. 
\qed 

\begin{COR}\label{lem:larger-component}
	Let $G$ be a $4.5$-tough  $(P_4\cup P_1)$-free  graph. Suppose that $C$ is a cycle of $G$ with order at least $\frac{3n}{t+1}$,
	and  $d_G(x) \ge \frac{3n}{t+1}$ for every vertex $x\in V(G)\setminus V(C)$. Then $G$ is Hamiltonian. 
\end{COR}

\pf We choose $C$ to be a longest cycle satisfying the conditions. If $C$ is Hamiltonian, then we are done. 
For otherwise, by Lemma~\ref{lem:vertex-insetting}, $G-V(C)$ has a component $H$ such that $|N_G(H, C)|<\frac{n}{t+1}$. 
Let $S=N_G(H, C)$. Then as  $d_G(x) \ge \frac{3n}{t+1}$ for every vertex $x\in V(G)\setminus V(C)$, it follows that 
$H$ is a component of $G-S$ of order at least $\frac{2n}{t+1}$.  Furthermore, as $|V(C)| \ge \frac{3n}{t+1}$ and $C-S$
is vertex-disjoint from $H$,  we know that   the  total number of vertices from components of  $G-S$ not containing 
a vertex of $H$ is at least $\frac{2n}{t+1}$.  Now,  by Theorem~\ref{thm:two-large-component-case}, $G$ is Hamiltonian. 
\qed

\section{Proof of Theorem~\ref{thm:main-result}}

We need the following result by 
H\"aggkvist and  Thomassen from 1982. 

\begin{THM}[{\cite[Theorem 2]{Cycles-through-edges}}]\label{thm:Hcycle-through-ind-edges}
Let $k\ge 0$ be an integer, and $G$ be a $(k+\alpha(G))$-connected graph, where $\alpha(G)$ is the independence number of $G$. 
Then for any linear forest $F$ of $G$ with at most $k$ edges, $G$ has a Hamiltonian cycle containing 
all the edges of $F$.  
\end{THM}

\proof[Proof of Theorem~\ref{thm:main-result}]  Let $n=|V(G)|$,  $S=\{v\in V(G): d_G(v) \ge \frac{n}{4}\}$,  and $T=V(G)\setminus S$.

\begin{CLA}\label{claim:no-p4}
	The graph $G-S$ is $P_4$-free. 
\end{CLA}

\pf  Assume otherwise that $G-S$ has an induced $P_4=u_1u_2u_3u_4$. Then as $G$ is $(P_4\cup P_1)$-free, it follows that 
$\max\{d_G(u_i): i\in [1,4]\} \ge \frac{n-4}{4}+1=\frac{n}{4}$, a contradiction to $u_i\not\in S$ for any $i$. 
\qed

Let $t=23$. We  may assume that  $G$ is not a complete graph. Thus $\delta(G) \ge 2t$ and so $n\ge 2t+1$.  
We consider two cases in  completing  the proof.

\medskip 

{\bf \noindent Case 1:  $|T|\ge \frac{3n}{t+1}$.}

If $G[T]$ has a Hamiltonian cycle, then  we are done by Corollary~\ref{lem:larger-component}.  
Thus we assume that $G[T]$ does not have  a Hamiltonian cycle. This, in particular, 
implies that $\delta(G[T]) <\frac{1}{2}|T|$ by Dirac's Theorem on Hamiltonian cycles.
Let $U\subseteq V(G[T])$ be a  minimum   cutset of $G[T]$. 
Then we have $|U| < \frac{1}{2}|T|$ and so $d_G(u, T\setminus U) =|T\setminus U|>\frac{1.5n}{t+1}$ for any $u\in U$ by Lemma~\ref{lem:P4-cut}(1). 
By Lemma~\ref{lem:S-vertex-large-degree-in-G-S}, we can find in $G$ a cycle $C$
containing all vertices of $T$ (an arbitrary ordering of vertices of $U$ plays the role of the ``ordering'' as specified in Lemma~\ref{lem:S-vertex-large-degree-in-G-S}).  Since $|V(C)| \ge \frac{3n}{t+1}$
and  all vertices of $G-V(C)$ have degree at least $\frac{n}{4}>\frac{3n}{t+1}$ in $G$, 
 Corollary~\ref{lem:larger-component} gives  a Hamiltonian cycle in $G$. 

\medskip 

{\bf \noindent Case 2:  $|T|<\frac{3n}{t+1}$.}

By Lemma~\ref{lem:path-partner-system1},  if $s(G-S) \le 0$, 
we find an $S$-matched   path-cover $\mathcal{Q}$ of $G-S$ with one single component; and if $s(G-S) \ge 1$, then 
we find an $S$-matched basic path-cover $\mathcal{Q}$ of $G-S$
with $ s(G-S)$ components.   As $G$ is $t$-tough, when $s(G-S) \ge 1$, we know that $c(\mathcal{Q}) \le \frac{n}{t+1}$.  
Let $k=\max\{1, s(G-S)\}$, and $x_iQ_iy_i$, where $x_i,y_i\in S$ for each $i\in [1,k]$,  be the $k$ components of $\mathcal{Q}$. 

We let $H$ be the  graph obtained from $G[S]$ by adding edges  $x_iy_i$ for each $i\in [1,k]$ whenever $x_iy_i\not\in E(G)$. 
Since $G$ is $t$-tough and so $\alpha(G) \le \frac{n}{t+1}$, we have $\alpha(H)  \le \frac{n}{t+1}$ as any independent set of $H$ is also an independent set of $G$.  
Furthermore, we have  $\delta(H) \ge \frac{n}{4}-|T|>\frac{3n}{t+1}$ by the definition of $S$.

Suppose first that $\frac{n}{4}-|T|-k-\frac{n}{t+1} \ge \frac{2n}{t+1}$. Under this assumption, we claim that 
 $H$ is $(k+\alpha(H))$-connected. For otherwise, let $W\subseteq V(H)$ be a minimum cutset.  Then $|W|< k+\alpha(H) \le \frac{2n}{t+1}$, and so 
 each component of $H-W$ has at least $\frac{n}{4}-|T|-|W| \ge \frac{2n}{t+1}$ vertices. 
Let $S^*=T\cup W$. Then $S^*$ is a cutset of $G$ such that $G-S^*$ has at least two components 
that each has order at least $\frac{2n}{t+1}$.  Applying  Theorem~\ref{thm:two-large-component-case},  we conclude that $G$ is Hamiltonian. 
Thus  we may assume that $H$ is $(k+\alpha(H))$-connected. Applying Theorem~\ref{thm:Hcycle-through-ind-edges}, $H$ has a Hamiltonian cycle $C$ 
going through all the edges $x_1y_1, \ldots, x_ky_k$. For each $i\in [1,k]$, by replacing each edge $x_iy_i$ on $C$ with the path $x_iQ_iy_i$, 
we obtain a Hamiltonian cycle of $G$.

We assume next that  $\frac{n}{4}-|T|-k-\frac{n}{t+1}<\frac{2n}{t+1}$. This gives  $|T|+2k>\frac{3n}{t+1}+k>\frac{3n}{t+1}$. 
We claim that $H$ is $(k+1)$-connected.  For otherwise, let $W\subseteq V(H)$ be a minimum cutset.  Then $|W| \le  k\le  \frac{n}{t+1}$, and so 
 each component of $H-W$ has at least $\frac{n}{4}-|T|-|W|>\frac{2n}{t+1}$ vertices. 
Let $S^*=T\cup W$. Then $S^*$ is a cutset of $G$ such that $G-S^*$ has at least two components 
that each has order at least $\frac{2n}{t+1}$.  Applying  Theorem~\ref{thm:two-large-component-case},  we conclude that $G$ is Hamiltonian.

Thus  $H$ is $(k+1)$-connected. 
By Theorem~\ref{thm:cycle-through-ind-edges}, $H$ has a cycle $C$ 
going through all the edges $x_1y_1, \ldots, x_ky_k$.  For each $i\in [1,k]$, by replacing each edge $x_iy_i$ on $C$ with the path $x_iQ_iy_i$, we get a  cycle $C^*$  in $G$
such that all vertices of $x_iQ_iy_i$ are covered by $C^*$.  As all the  $k$ paths $x_1Q_1y_1, \ldots, x_kQ_ky_k$  together cover all the vertices of $T$ and $2k$ vertices from $S$, 
we know that the order of $C^*$ is at least $\frac{3n}{t+1}$. We also have   $V(G)\setminus V(C^*) \subseteq S$.  
Now we find in $G$ a  Hamiltonian cycle again by 
Corollary~\ref{lem:larger-component}. 
\qed

\section*{Acknowledgment}
The author would like to thank Masahiro Sanka for reading an early draft and providing some  feedback, and the two anonymous referees for their careful reading and valuable comments.
The author  was  partially supported by NSF grants  DMS-2345869 and DMS-2451895.

\bibliographystyle{abbrv}
\bibliography{SSL}

\end{document}